\documentclass[12pt,a4paper]{article}
\usepackage{amssymb,amsmath,amsthm,euscript}
\usepackage{leftidx}
\usepackage[matrix,arrow,curve]{xy}

\usepackage[utf8]{inputenc}
\usepackage[english]{babel}

\newcommand{\chara}{\mathop{\mathrm{char}}\nolimits}
\newcommand{\End}{\mathop{\mathrm{End}}\nolimits}
\newcommand{\Ker}{\mathop{\mathrm{Ker}}\nolimits}

\newcommand{\Hom}{\mathop{\mathrm{Hom}}\nolimits}
\newcommand{\bfhom}{\mathop{\mathrm{\bf hom}}\nolimits}
\newcommand{\bfend}{\mathop{\mathrm{\bf end}}\nolimits}
\newcommand{\bfcohom}{\mathop{\mathrm{\bf cohom}}\nolimits}
\newcommand{\bfcoend}{\mathop{\mathrm{\bf coend}}\nolimits}
\newcommand{\Mat}{\mathop{\mathrm{Mat}}\nolimits}

\newcommand{\id}{\mathop{\mathrm{id}}\nolimits}

\newcommand{\ootimes}{\mathbin{\underline{\otimes}}}
\newcommand{\isoright}{\xrightarrow{\smash{\raisebox{-0.65ex}{\ensuremath{\sim}}}}}
\newcommand{\xisoright}[1]{\xrightarrow[\smash{\raisebox{0.65ex}{\ensuremath{\sim}}}]{#1}}

\numberwithin{equation}{section}
\renewcommand{\Im}{\mathop{\mathrm{Im}}\nolimits}

\newcommand{\op}{{\mathrm{op}}}
\newcommand{\cop}{{\mathrm{cop}}}

\renewcommand{\le}{\leqslant}
\renewcommand{\ge}{\geqslant}

\newcommand{\A}{{\mathcal A}}
\newcommand{\B}{{\mathcal B}}
\newcommand{\C}{{\mathcal C}}

\newcommand{\M}{{\mathcal M}}

\newcommand{\R}{{\mathcal R}}
\newcommand{\U}{{\mathcal U}}

\newcommand{\gX}{{\mathfrak X}}

\newcommand{\bz}{{\bar0}}
\newcommand{\bu}{{\bar1}}
\newcommand{\wt}{\widetilde}
\newcommand{\wh}{\widehat}
\newcommand{\st}{{\mathbf{st}}}
\newcommand{\ist}{{\mathbf{ist}}}

\newtheorem{Th}{Theorem}[section]
\newtheorem{Lem}[Th]{Lemma}
\newtheorem{Prop}[Th]{Proposition}

\newtheorem{Def}[Th]{Definition}

\newcommand{\BB}{{\mathbb{B}}}

\newcommand{\ZZ}{{\mathbb{Z}}}
\newcommand{\MM}{{\mathbb{M}}}
\newcommand{\OO}{{\mathbb{O}}}
\newcommand{\NN}{{\mathbb{N}}}
\newcommand{\KK}{{\mathbb{K}}}

\renewcommand{\AA}{{\mathbb{A}}}

\newcommand{\Set}{{\mathbf{Set}}}

\newcommand{\AffSch}{{\mathbf{AffSch}}}
\newcommand{\SAffSch}{{\mathbf{SAffSch}}}

\newcommand{\Vect}{{\mathbf{Vect}}}
\newcommand{\FVect}{{\mathbf{FVect}}}
\newcommand{\Alg}{{\mathbf{Alg}}}
\newcommand{\CommAlg}{{\mathbf{CommAlg}}}
\newcommand{\FAlg}{{\mathbf{FAlg}}}
\newcommand{\SVect}{{\mathbf{SVect}}}
\newcommand{\SAlg}{{\mathbf{SAlg}}}
\newcommand{\CommSAlg}{{\mathbf{CommSAlg}}}
\newcommand{\FSVect}{{\mathbf{FSVect}}}
\newcommand{\FSAlg}{{\mathbf{FSAlg}}}
\newcommand{\GrVect}[1]{{\text{${#1}$-$\mathbf{GrVect}$}}}
\newcommand{\GrAlg}[1]{{\text{${#1}$-$\mathbf{GrAlg}$}}}
\newcommand{\GrSVect}[1]{{\text{${#1}$-$\mathbf{GrSVect}$}}}
\newcommand{\GrSAlg}[1]{{\text{${#1}$-$\mathbf{GrSAlg}$}}}

\newcommand{\QA}{{\mathbf{QA}}}
\newcommand{\FQA}{{\mathbf{FQA}}}
\newcommand{\QSA}{{\mathbf{QSA}}}
\newcommand{\FQSA}{{\mathbf{FQSA}}}
\newcommand{\CommQSA}{{\mathbf{CommQSA}}}

\newcommand{\bfC}{{\mathbf{C}}}
\newcommand{\bfD}{{\mathbf{D}}}

\newcommand{\bfP}{{\mathbf{P}}}

\newcommand{\bfk}{{\mathbf{k}}}
\newcommand{\bfl}{{\mathbf{l}}}

\newcommand{\Mon}{\mathop{\mathbf{Mon}}\nolimits}
\newcommand{\cMon}{\mathop{\mathbf{\leftidx{_c}Mon}}\nolimits}
\newcommand{\Comon}{\mathop{\mathbf{Comon}}\nolimits}
\newcommand{\cocComon}{\mathop{\mathbf{\leftidx{_{coc}}Comon}}\nolimits}
\newcommand{\Bimon}{\mathop{\mathbf{Bimon}}\nolimits}
\newcommand{\Lact}{\mathop{\mathbf{Lact}}\nolimits}
\newcommand{\Rep}{\mathop{\mathbf{Rep}}\nolimits}
\newcommand{\Corep}{\mathop{\mathbf{Corep}}\nolimits}
\newcommand{\Lcoact}{\mathop{\mathbf{Lcoact}}\nolimits}

\usepackage{hyperref}

\newcounter{bbcount}[subsection]
\renewcommand{\thebbcount}{{\thesection}.\arabic{subsection}.\arabic{bbcount}}

\newcommand{\bb}[1]{\vspace{2mm}\noindent\refstepcounter{bbcount}{\bf\thebbcount.}{ \bfseries #1}}

\title{Quantum Representation Theory and \\ Manin matrices~II: super case}
\author{Alexey Silantyev\thanks{aleksejsilantjev@gmail.com}}
\date{}

\begin{document}

\maketitle

\vspace{-5mm}
\begin{center}
{\it Bogoliubov Laboratory of Theoretical Physics, Joint Institute for Nuclear Research, 141980~Dubna, Moscow region, Russia} \\
{\it State University ``Dubna''{}, Universitetskaya~st.~19, 141980~Dubna, Moscow region, Russia} \\
\end{center}
\vspace{5mm}

\begin{abstract}
 We construct super-version of Quantum Representation Theory. The quadratic super-algebras and operations on them are described. We also describe some important monoidal functors. We proved that the monoidal category of graded super-algebras with Manin product is coclosed relative to the subcategory of finitely generated quadratic super-algebras. The super-version of the $(A,B)$-Manin matrices is introduced and related with the quadratic super-algebras. We define a super-version of quantum representations and of quantum linear actions, relate them to each other and describe them by the super-Manin matrices. Some operations on quantum representations/quantum linear actions are described. We show how the classical representations lift to the quantum level.
\end{abstract}



\tableofcontents

\section{Introduction}

It is well-known that in the usual differential/algebraic geometry the spaces contravariantly corresponds to the algebras of functions. The idea to generalise different kinds of spaces by substituting the commuting algebras of functions by non-commutative algebras is quite old. One source of such generalisation is quantisation in physics. By this reason such `non-commutative' generalisation of a space are often called `quantum space'. A generalisation of Lie/algebraic group was called `quantum group' by Drinfeld~\cite{Dr}.

In~\cite{Manin88} Yuri Manin proposed to consider quadratic algebras as quantum analogues of the vector spaces. This was a starting point to generalise classical representation theory to the quantum case in~\cite{Sqrt}; namely, the notion of representation was generalised to the case of a quantum representation space. Such quantum representations are described by so-called Manin matrices.

In some form the Manin matrices appeared in~\cite{Manin87,Manin88} in connection with the quadratic algebras. The Manin matrices for polynomial algebras and their $q$-deformations were investigated and applied in the works~\cite{GLZ,CF,CM,CFR,RST,qManin,IO,MolevSO}. For general quadratic algebras they were described in~\cite{Smm}. Super-analogues of the Manin matrices was considered in~\cite{Manin89,MR}.

The main tool of Quantum Representation Theory in~\cite{Sqrt} is the theory of monoidal categories. The correspondence between linear actions and representations in a closed monoidal category was generalised to a wider class of monoidal categories. This general representation theory was developed in~\cite{Sgrt}.

In the present paper we consider super-versions of quadratic algebras, of Manin matrices and of operations with the quadratic algebras. We apply the general approach of~\cite{Sgrt} to an appropriate monoidal category and describe quantum representations by the super-Manin matrices.

The article is organised as follows. Section~\ref{sec2} is preliminary: we recall and introduce some terms and notations here. The quadratic super-algebras are investigated in Section~\ref{sec3}. The binary operations (monoidal products) on quadratic super-algebras are described in Subsection~\ref{sec31}. In Subsection~\ref{sec32} we collect useful functors and describe their properties. Subsection~\ref{sec33} is devoted to the coclosed structure on a monoidal category of the quadratic super-algebras given by so-called internal cohom-functor. In Subsection~\ref{sec41} we introduce super-Manin matrices for a pair of quadratic super-algebras in terms idempotent operators; also we introduce universal super-Manin matrices and use them to describe the internal cohom-functor. Quantum Representation Theory for the super-case is described in Subsection~\ref{sec42}.

\vspace{3mm}
{\it Acknowledgements}. 
The author is grateful to A.~Isaev, Yu.~Manin and E.~Patrin for references. The author also thanks V.~Rubtsov for useful advice.

\section{Preliminaries}
\label{sec2}

In this work we use the terms and notations introduced in~\cite[\S~2]{Sqrt} and $\cite[\S~2]{Sgrt}$. Fist we briefly remind some of them. Then we introduce the structure of $\AA$-grading and consider the main category we use -- the category of $\AA$-graded algebras, where $\AA=\ZZ_2\times\ZZ$.

\subsection{Notations and conventions}
\label{sec21}

\bb{Monoidal categories.} We denote a monoidal category by $(\bfC,\otimes)$ or $(\bfC,\otimes,I)$, where $\bfC$ is a category with a bifunctor $\otimes\colon\bfC\times\bfC\to\bfC$ and a unit object $I=I_\bfC$. We suppose monoidal categories to be strict. The monoids and comonoids in this monoidal category form the categories $\Mon(\bfC,\otimes)$ and $\Comon(\bfC,\otimes)=\big(\Mon(\bfC^\op,\otimes)\big)^\op$. For any $\MM\in\Mon(\bfC,\otimes)$ and $\OO\in\Comon(\bfC,\otimes)$ we have categories of (left) actions and coactions $\Lact(\MM)$ and $\Lcoact(\OO)$.

Let the monoidal category $(\bfC,\otimes)$ be symmetric. Then $\Mon(\bfC,\otimes)$ and $\Comon(\bfC,\otimes)$ are monoidal and also symmetric (the monoidal product in this category is usually denoted by the same symbol). They have full subcategories of commutative monoids and cocommutative comonoids: $\cMon(\bfC,\otimes)\subset\Mon(\bfC,\otimes)$, $\cocComon(\bfC,\otimes)\subset\Mon(\bfC,\otimes)$. In this case one can also define bimonoids, which form the category
$$\Bimon(\bfC,\otimes)=\Comon\big(\Mon(\bfC,\otimes),\otimes\big)=\Mon\big(\Comon(\bfC,\otimes),\otimes\big).$$
Any bimonoid $\BB\in\Bimon(\bfC,\otimes)$ gives the monoidal categories $\Lact(\BB)$ and $\Lcoact(\BB)$. They are symmetric if $\BB$ is cocommutative or commutative respectively.

Any category $\bfC$ with finite products is a symmetric monoidal category $(\bfC,\times)$, where $\times$ is the categorical product and the unit object is the terminal object. An example of such category is the category of sets $\Set$.

\bb{Monoidal functors.} A {\it lax} monoidal structure of a functor $F\colon\bfC\to\bfD$ between two monoidal categories $(\bfC,\otimes)$ and $(\bfD,\odot)$ is given by morphisms $\varphi\colon I_\bfD\to FI_\bfC$ and $\phi_{X,Y}\colon FX\odot FY\to F(X\otimes Y)$ satisfying some conditions. A {\it colax} monoidal structure of $F\colon\bfC\to\bfD$ is given by $\varphi\colon FI_\bfC\to I_\bfD$ and $\phi_{X,Y}\colon F(X\otimes Y)\to FX\odot FY$. Lax/colax monoidal functor is the triple $F=(F,\varphi,\phi)\colon(\bfC,\otimes)\to(\bfD,\odot)$. A monoidal functor between symmetric monoidal categories is called {\it symmetric} if it respects the symmetric structures.

Any (co)lax monoidal functor translates (co)monoids to (co)monoids and (co)actions to (co)actions, so we have the induced functors $\Mon(F)\colon\Mon(\bfC,\otimes)\to\Mon(\bfD,\odot)$, $F_\MM\colon\Lact(\MM)\to\Lact\big(\Mon(F)\MM\big)$ for the lax case and  $\Comon(F)\colon\Comon(\bfC,\otimes)\to\Comon(\bfD,\odot)$, $F^\OO\colon\Lcoact(\OO)\to\Lcoact\big(\Comon(F)\OO\big)$ for the colax case.

If all $\varphi$ and $\phi_{X,Y}$ are isomorphisms, then the monoidal functor $F=(F,\varphi,\phi)\colon(\bfC,\otimes)\to(\bfD,\odot)$ is called {\it strong} monoidal. A symmetric strong monoidal functor induces the functors $\Bimon(F)\colon\Mon(\bfC,\otimes)\to\Bimon(\bfD,\odot)$, $F_\BB\colon\Lact(\BB)\to\Lact\big(\Mon(F)\BB\big)$ and $F^\BB\colon\Lcoact(\BB)\to\Lcoact\big(\Comon(F)\BB\big)$.

A contravariant functor $F\colon\bfC\to\bfD$ is called (symmetric) lax/colax/strong monoidal iff the corresponding covariant functor $\bar F\colon\bfC^\op\to\bfD$ is (symmetric) lax/colax/strong monoidal (or, equivalently, the opposite functor $\bar F^\op\colon\bfC\to\bfD^\op$ is (symmetric) colax/lax/strong monoidal).

\bb{Relative adjoints}~\cite{Ulm}, \cite[\S~2.3]{Sgrt}. Consider categories $\bfC$ and $\bfD$. We say that a functor $F\colon\bfC\to\bfD$ has a {\it right adjoint relative} to a full subcategory $\bfD'\subset\bfD$ iff there exists a functor $G\colon\bfD'\to\bfC$ (called right adjoint for $F$ relative to $\bfD'$) and a bijection
\begin{align} \label{HomFXGZ}
 \Hom(FX,Z)\cong\Hom(X,GZ)
\end{align}
natural in $X\in\bfC$, $Z\in\bfD'$. We say that a functor $G\colon\bfD\to\bfC$ has a {\it left adjoint relative} to a full subcategory $\bfC'\subset\bfC$ iff there exists a functor $F\colon\bfC'\to\bfD$ (called left adjoint for $G$ relative to $\bfC'$) and a bijection~\eqref{HomFXGZ} natural in $X\in\bfC'$ and $Z\in\bfD$. The right/left adjoint with the {\it relative adjunction}~\eqref{HomFXGZ} is unique up to an isomorphism. A functor $F\colon\bfC\to\bfD$ has a right adjoint relative to $\bfD'\subset\bfD$ iff there exist a universal morphism from $F$ to each object $Z\in\bfD'$. A functor $G\colon\bfD\to\bfC$ has a left adjoint relative to $\bfC'\subset\bfC$ iff there exist universal morphism from each object $X\in\bfC'$ to $G$.

\bb{Relatively closed monoidal categories.} In~\cite[\S~2.4]{Sgrt} we generalised the notions of closed and coclosed monoidal category in terms of relative adjoints. We call a monoidal category $(\bfC,\otimes)$ {\it closed/coclosed relative} to a full subcategory $\bfP\subset\bfC$ or {\it relatively closed/coclosed with parametrising subcategory} $\bfP$ iff for any $Y\in\bfP$ the functor $F_Y=-\otimes Y$ has a right/left adjoint relative to $\bfP$.

Fix parametrising subcategory $\bfP$. Then $(\bfC,\otimes)$ is relatively closed iff there exists a bifunctor $\bfhom\colon\bfP^\op\times\bfP\to\bfC$ (called internal hom-functor) and a bijection
\begin{align}
 \theta\colon\Hom\big(X,\bfhom(Y,Z)\big)=\Hom(X\otimes Y,Z)
\end{align}
natural in $X\in\bfC$, $Y,Z\in\bfP$. In this case the internal end-object $\bfend(Y)=\bfhom(Y,Y)$ has a structure of monoid in $(\bfC,\otimes)$ for any $Y\in\bfP$.

Dually, the monoidal category $(\bfC,\otimes)$ is relatively coclosed iff there exists a bifunctor $\bfcohom\colon\bfP^\op\times\bfP\to\bfC$ (called internal cohom-functor) and a bijection
\begin{align}
 \vartheta\colon\Hom\big(\bfcohom(Y,X),Z\big)=\Hom(X,Z\otimes Y)
\end{align}
natural in $X,Y\in\bfP$, $Z\in\bfC$. In this case the internal coend-object $\bfcoend(Y)=\bfcohom(Y,Y)$ has a structure of comonoid in $(\bfC,\otimes)$ for any $Y\in\bfP$.

The monoidal category is closed/coclosed iff it is closed/coclosed relative to the whole category $\bfC$.

\bb{Vector spaces and algebras.} Fix an infinite field $\KK$ of characteristic $\chara\KK\ne2$. We consider all vector spaces over $\KK$. The category of such vector spaces is denoted by $\Vect$. This is a monoidal category with the standard tensor product $V\otimes W=V\otimes_\KK W$. The monoids in this category form the category $\Alg=\Mon(\Vect,\otimes)$. We call its objects simply {\it algebras} (these are associative unital algebras over $\KK$). An identity element of an algebra $\A$ is denoted by $1_\A$ or $1$.

The category $\Alg$ is also monoidal with respect to the tensor product of algebras. Its full monoidal subcategory consisting of the commutative algebras is opposite to the category of the affine schemes over $\KK$: $(\CommAlg,\otimes)^\op=(\AffSch,\times)$. Denote the categories of finite-dimensional vector spaces and algebras by $\FVect$ and $\FAlg=\Mon(\FVect,\otimes)$. These are monoidal subcategories of $(\Vect,\otimes)$ and $(\Alg,\otimes)$.

\bb{$\AA$-graded vector spaces.} Let $\AA\in\cMon(\Set,\times)$ be an Abelian monoid with the binary operation $\AA\times\AA\to\AA$ written additively: $(a,b)\mapsto a+b$. Denote by $\GrVect{\AA}$ the category of $\AA$-graded vector spaces $V=\bigoplus\limits_{g\in\AA}V_g$. Elements of the component $V_g$ are called homogeneous of degree $g$. The morphisms in this category are graded linear maps $f\colon V\to W$, i.e. $f(V_g)\subset W_g$. Each such map is given by a family of arbitrary linear maps $f_g\colon V_g\to W_g$. In other words, $\GrVect{\AA}$ is the category of functors $\AA\to\Vect$.

The tensor product of two $\AA$-graded vector spaces decomposed as
\begin{align}
 &V\otimes W=\bigoplus_{g\in\AA}(V\otimes W)_g, &&(V\otimes W)_g=\bigoplus_{h,h'\in \AA\atop h+h'=g}V_{h}\otimes W_{h'}
\end{align}
(if $\AA$ is a group then $(V\otimes W)_g=\bigoplus\limits_{h\in\AA}V_{h}\otimes W_{g-h}$). In this way we obtain a monoidal category $(\GrVect{\AA},\otimes)$. Its unit object is the vector space $\KK$ with the grading $(\KK)_0=\KK$, $\KK_g=0$ for $g\ne 0$, where $0$ is the neutral element of $\AA$.

The monoidal category $(\GrVect{\AA},\otimes)$ is closed. The component $\bfhom(W,Z)_{g}$ is the vector space consisting of the operators $f\colon W\to Z$ such that $f(W_h)=V_{g+h}$.

\bb{$\AA$-graded algebras.} A monoid $\A\in\Mon(\GrVect{\AA},\otimes)$ is an algebra with a grading $\A=\bigoplus\limits_{g\in\AA}\A_g$ such that $A_gA_h\subset A_{g+h}$ and $1_\A\in\A_0$ (if $\AA$ is a group or a submonoid of a group, then the former condition implies the latter one). We obtain the category $\GrAlg{\AA}=\Mon(\GrVect{\AA},\otimes)$, where the morphisms are graded homomorphisms of algebras. Its objects are called $\AA$-graded algebras.

If $(\GrVect{\AA},\otimes)$ is equipped with a symmetric structure, then it gives a monoidal structure on $\GrAlg{\AA}$ with the symmetric structure.
Recall that the standard symmetric structure of $(\Vect,\otimes)$ is defined by the permutation operators $\sigma_{V,W}\colon V\otimes W\to W\otimes V$, $v\otimes w\mapsto v\otimes w$. These operators are graded due to the commutativity of $\AA$, so it gives a simplest symmetric structure on $(\GrVect{\AA},\otimes)$. However, it is not unique.

Let $\AA'\subset\AA$ be a submonoid. Then $\GrVect{\AA'}$ and $\GrAlg{\AA'}$ are full monoidal subcategories of $\GrVect{\AA}$ and $\GrAlg{\AA}$, consisting of the $\AA$-graded vector spaces (algebras) $V$ such that $V_g=0$ for $g\in\AA\backslash\AA'$. In particular, $\Vect=\GrVect{\{0\}}\subset\GrVect{\AA}$, $\Alg\subset\GrAlg{\AA}$.

\subsection{Graded super-algebras}
\label{sec22}

\bb{Super-vector spaces and super-algebras.} \label{bbSVect}
Let $\ZZ_2=\ZZ/2\ZZ=\{\bz,\bu\}$. This is an Abelian additive group, which has a structure of ring given by the multiplication $\bar k\cdot\bar l=\overline{kl}$. Consider the category $\SVect:=\GrVect{\ZZ_2}$. Its objects $V\in\SVect$ are called {\it super-vector spaces}. The homogeneous elements $v\in V_\bz$ are called {\it even}, while $v\in V_\bu$ are called {\it odd}. For $v\in V_{\bar k}$ we use the notation $[v]:=\bar k$.

A subspace $\wt V$ of a super-vector space $V\in\SVect$ is called super-subspace iff $\wt V$ has a structure of $\ZZ_2$-grading such that the embedding $\wt V\hookrightarrow V$ is $\ZZ_2$-graded, i.e. $\wt V=\wt V_\bz\oplus\wt V_\bu$ where $\wt V_\bz=\wt V\cap V_\bz$ and $\wt V_\bu=\wt V\cap V_\bu$.

The objects of $\SAlg:=\GrAlg{\ZZ_2}$ are monoids in the monoidal category $(\SVect,\otimes)$ called {\it super-algebras}. The super-vector space $\bfend(V)=\bfhom(V,V)$ is equipped with a structure of super-algebra.
 The categories of finite-dimensional super-vector spaces and super-algebras are full monoidal subcategories $\FSVect$ and $\FSAlg=\Mon(\FVect,\otimes)$ of $(\SVect,\otimes)$ and $(\SAlg,\otimes)$ respectively.

Fix the following symmetric structure on this monoidal category:
\begin{align} \label{sigmaSuper}
 \sigma(v\otimes w)=(-1)^{[v][w]}w\otimes v,
\end{align}
where $(-1)^{\bar k}=(-1)^k$ (the usage of the notation $[v]$ supposes that $v$ is homogeneous). We obtain the symmetric monoidal category $(\GrAlg{\ZZ_2},\otimes)$. The product $\R\otimes\wt\R$ of two super-algebras $\R$ and $\wt\R$ is defined by the formula
\begin{align}
 &(a\otimes b)\cdot(c\otimes d)=(-1)^{[b][c]}(a\cdot c)\otimes(b\cdot d), &&a,c\in\R,\quad b,d\in\wt\R.
\end{align}
The unit object of $(\GrAlg{\ZZ_2},\otimes)$ is the algebra $\KK\in\Alg\subset\SAlg$, the symmetric structure is defined by the same formula~\eqref{sigmaSuper}. A commutative super-algebra is a super-algebra $\R$ satisfying the condition $ab=(-1)^{[a][b]}ba$, $a,b\in\R$. The commutative super-algebras form a full monoidal subcategory $\CommSAlg\subset\SAlg$, where the functor $\otimes$ coincides with categorical coproduct, so $(\CommSAlg,\otimes)$ is opposite to $(\SAffSch,\times)$, where $\SAffSch=\CommSAlg^\op$ is the category of the affine super-schemes over $\KK$. The objects of $\SAlg^\op$ are the quantum analogues of the affine super-schemes.

Affine algebraic super-group is a group in $\SAffSch$. More generally affine algebraic super-monoid is a monoid in $\SAffSch$, i.e. a monoid in $(\CommSAlg^\op,\otimes)$. The monoids in $(\SAlg^\op,\otimes)$ are called {\it quantum super-monoids}. Their category is opposite to the category of super-bialgebras: $\Mon(\SAlg^\op,\otimes)=\Bimon(\SVect^\op,\otimes)=\Bimon(\SVect,\otimes)^\op$. The quantum super-groups correspond to the Hopf super-algebras.

The symmetric monoidal category $(\SVect,\otimes)$ is closed. For $V,W\in\SVect$ the internal hom-object $\bfhom(V,W)$ is a super-vector space consisting of all the linear maps $f\colon V\to W$. The elements $f\in\bfhom(V,W)_\bz$ preserve the $\ZZ_2$-grading and we identify them with the morphisms $f\in\Hom(V,W)$. The elements $f\in\bfhom(V,W)_\bu$ are operators $f\colon V\to W$ such that $f(V_\bz)\subset W_\bu$ and $f(V_\bu)\subset W_\bz$.

Note that direct sum of vector spaces provides another monoidal product on the category $\SVect$. The direct sum of $V,W\in\SVect$ is a super-vector space $V\oplus W$ with the components $(V\oplus W)_\bz=V_\bz\oplus W_\bz$, $(V\oplus W)_\bu=V_\bu\oplus W_\bu$. The monoidal category $(\SVect,\oplus)$ has a canonical symmetric structure. 

\bb{$\AA$-graded super-algebras.} Consider the product of an Abelian monoid $\AA$ and the group $\ZZ_2$. The $\AA\times\ZZ_2$-grading provides the categories $\GrSVect{\AA}:=\GrVect{\AA\times\ZZ_2}$ and $\GrSAlg{\AA}:=\GrAlg{\AA\times\ZZ_2}$; we call their objects {\it $\AA$-graded super-vector spaces} and {\it $\AA$-graded super-algebras} respectively. Each object $V\in\GrSVect{\AA}$ has decompositions
\begin{align}
 V=\bigoplus_{g\in\AA}(V_{g\bz}\oplus V_{g\bu})=\bigoplus_{g\in\AA}V_g=V_\bz\oplus V_\bu,
\end{align}
where $V_{g\bar k}=V_{(g,\bar k)}\in\Vect$, $V_g=V_{g\bz}\oplus V_{g\bu}\in\GrVect{\ZZ_2}$, $V_{\bar k}=\bigoplus\limits_{g\in\AA}V_{g\bar k}\in\GrVect{\AA}$. Thus an $\AA$-graded super-vector space (super-algebra) is a vector space (an algebra) equipped with two gradings ($\AA$- and $\ZZ_2$-grading) which are compatible with each other. We have the following commutative diagrams of the forgetful functors:
\begin{align} \label{ASVectForg}
&\xymatrix{
 \GrSVect{\AA}\ar[r]\ar[d] &\SVect\ar[d] \\
 \GrVect{\AA}\ar[r] &\Vect
} &
&\xymatrix{
 \GrSAlg{\AA}\ar[r]\ar[d] &\SAlg\ar[d] \\
 \GrAlg{\AA}\ar[r] &\Alg
}
\end{align}
Note that $\Vect=\GrVect{\{0\}}$. The embedding of the monoids gives the commutative diagrams of the embedding functors
\begin{align} \label{ASVectEmb}
&\xymatrix{
 \Vect\ar@{^{(}->}[r]\ar@{^{(}->}[d] &\SVect\ar@{^{(}->}[d] \\
 \GrVect{\AA}\ar@{^{(}->}[r] &\GrSVect{\AA}
} &
&\xymatrix{
 \Alg\ar@{^{(}->}[r]\ar@{^{(}->}[d] &\SAlg\ar@{^{(}->}[d] \\
 \GrAlg{\AA}\ar@{^{(}->}[r] &\GrSAlg{\AA}
}
\end{align}

\bb{Monoidal products of $\AA$-graded super-algebras.}
Note that $(\GrSVect{\AA},\otimes)$ and $(\GrSAlg{\AA},\otimes)=\Mon(\GrSVect{\AA},\otimes)$ are symmetric monoidal categories with the standard tensor product of vector spaces and symmetric structure given by~\eqref{sigmaSuper}. The functors in~\eqref{ASVectEmb} are symmetric strong monoidal with respect to $\otimes$. The functors in the left diagram~\eqref{ASVectForg} are also strong monoidal, but the vertical ones are not symmetric. By this reason the vertical functors in the right diagram~\eqref{ASVectForg} are not monoidal.

Let us introduce another monoidal structure on $\GrSAlg{\AA}$: for two $\AA$-graded super-algebras $\A$ and $\B$ we define {\it Manin (white) product} by the formula
\begin{align}
 &\A\circ\B=\bigoplus_{g\in\AA}(\A\circ\B)_g, &&(\A\circ\B)_g=\A_g\otimes\B_g;
\end{align}
as a super-graded vector space it has the components $(\A\circ\B)_\bz=\bigoplus\limits_{g\in\AA}(\A_{g\bz}\otimes\B_{g\bz}\oplus\A_{g\bu}\otimes\B_{g\bu})$, $(\A\circ\B)_\bu=\bigoplus\limits_{g\in\AA}(\A_{g\bz}\otimes\B_{g\bu}\oplus\A_{g\bu}\otimes\B_{g\bz})$.
Note that $\A\circ\B$ is a super-subalgebra in $\A\otimes\B$, but the embedding is not graded. We also use the notation $V\circ W$ for $\AA$-graded super-vector space with components $V_g\circ W_g$, where $V,W\in\GrSVect{\AA}$.

The cases we are interested in are $\AA=\NN_0$ and $\AA=\ZZ$, where $\NN_0$ is a submonoid of $\ZZ$, consisting of non-negative integers. The Manin product for $\NN_0$- and $\ZZ$-graded algebras is defined in~\cite{Sqrt}, \cite{Sgrt}. The embedding $\GrAlg{\NN_0}\hookrightarrow\GrSAlg{\NN_0}$ is symmetric strong monoidal with respect to the Manin product $\circ$. By this reason the super-version of Quantum Representation Theory we present here includes the corresponding results of~\cite{Sqrt}.

Note that the unit objects of $(\GrSAlg{\NN_0},\circ)$ and $(\GrSAlg{\ZZ},\circ)$ are the graded algebras of polynomials $\KK[u]\in\GrAlg{\NN_0}$ and of Laurent polynomials $\KK[u,u^{-1}]\in\GrAlg{\ZZ}$ respectively. The embedding functor $\GrSAlg{\NN_0}\hookrightarrow\GrSAlg{\ZZ}$ is only colax monoidal: $\phi_{\A,\B}=\id_{\A\otimes\B}$ are isomorphisms, but $\varphi\colon\KK[u]\hookrightarrow\KK[u,u^{-1}]$ is not isomorphism.

\section{Categories of quadratic super-algebras}
\label{sec3}

The main ingredient of Quantum Representation Theory is the category of representation spaces. In the super-case this is the category of the (finitely generated) quadratic super-algebras (more precisely, its opposite category). Here we describe different types of operations on quadratic super-algebras: binary operations (monoidal products), non-binary operations realised as functors from the category of quadratic algebras. We also derive the internal (co)hom-functors, which play a key role in the theory.

\subsection{Monoidal structures}
\label{sec31}

We use the standard tensor notations $f^{(a)}$ for a morphism $f$ putted to the $a$-th tensor (monoidal) factor. If $f$ acts on a monoidal product of two objects, then we use the notation $f^{(a,a+1)}$. For example, $$\sigma^{(23)}=\id\otimes\sigma\otimes\id\otimes\cdots\otimes\id\colon V_1\otimes V_2\otimes V_3\otimes V_4\otimes\cdots\otimes V_k\to V_1\otimes V_3\otimes V_2\otimes V_4\otimes\cdots\otimes V_k,$$ where $\sigma=\sigma_{V_2,V_3}\colon V_2\otimes V_3\to V_3\otimes V_2$ is the symmetric structure.
We need the following fact from the linear algebra (see e.g.~\cite[\S~1, ex.~1.3]{LQ}).

\begin{Lem} \label{LemTPQ}
 For any vector spaces $V,W\in\Vect$ and subspaces $V_0\subset V$ and $W_0\subset W$ there is an isomorphism
\begin{align} \label{VV0WW0}
 (V/V_0)\otimes(W/W_0)\isoright(V\otimes W)/(V_0\otimes W+V\otimes W_0)
\end{align}
which sends $\bar v\otimes\bar w$ to the class $\overline{v\otimes w}$, where $\bar v\in V/V_0$ is the class of $v\in V$.
If $V,W$ have $\ZZ_2$-gradings and $V_0$, $W_0$ are super-subspaces, then the isomorphism~\eqref{VV0WW0} is $\ZZ_2$-graded.
\end{Lem}

\noindent{\bf Proof.} It is obvious that $\bar v\otimes\bar w\mapsto\overline{v\otimes w}$ is a correctly defined epimorphism of super-vector spaces. By choosing appropriate bases in $V$ and $W$ one can show that its kernel vanishes. \qed

\bb{Quadratic super-algebras.}
Let $V\in\SVect$.  The tensor algebra $TV$ has a structure of super-algebra compatible with the $\NN_0$-grading: $TV\in\GrSAlg{\NN_0}$. Its components with respect to the $\NN_0$-grading are $(TV)_0=\KK$, $(TV)_1=V$, $(TV)_2=V\otimes V$, $(TV)_k=V^{\otimes k}$, $k\in\NN_0$.

The quotient of $TV$ over an ideal $(R)$ generated by a super-subspace $R\subset V\otimes V$ is also an $\NN_0$-graded super-algebra. A {\it quadratic super-algebra}
\footnote{More generally, one can define quadratic super-algebras over an algebra $\R$ as $\NN_0$-graded super-algebras $\A$ generated by $\A_1$ over $\A_0=\R$. By default we mean the quadratic super-algebras over $\KK$ (connected quadratic super-algebras), cf.~\cite{Smm}, \cite{Sqrt}.
}
 is an $\NN_0$-graded super-algebra $\A$ generated by elements of $\A_1$ with quadratic commutation relations, that is $\A\cong TV/(R)$ for some $V\in\SVect$ and super-subspace $R\subset V\otimes V$. The category of quadratic super-algebras form a full subcategory $\QSA\subset\GrSAlg{\NN_0}$. The quadratic super-algebra $TV/(R)$ is finitely generated iff $V\in\FSVect$. Denote the full subcategory of the finitely generated quadratic super-algebras by $\FQSA\subset\QSA$.

\bb{Super-version of Manin's binary operations.}
 Consider $\QSA$ as a full subcategory of $\GrSAlg{\NN_0}$. The following fact generalises~\cite[\SS~3, Lemma~7]{Manin88} for the super-case.

\begin{Prop} \label{Propotimescirc}
 The subcategory $\QSA\subset\GrSAlg{\NN_0}$ is monoidal with respect to the both monoidal products $\otimes$ and $\circ$. Let $\A=TV/(R)$ and $\B=TW/(S)$, where $V,W\in\SVect$ and $R\subset V^{\otimes2}$, $S\subset W^{\otimes2}$ are super-subspaces. Then we have isomorphisms
\begin{align}
 &\A\otimes\B\cong T(V\oplus W)/(R'), &&R'=R\oplus S\oplus[V,W], \label{AotimesB} \\
 &\A\circ\B\cong T(V\otimes W)/(R_{\mathrm w}), &&R_{\mathrm w}=\sigma^{(23)}(R\otimes W^{\otimes2}+V^{\otimes2}\otimes S), \label{AcircB}
\end{align}
where $[V,W]\subset V\otimes W\oplus W\otimes V\subset(V\oplus W)^{\otimes2}$ is spanned by $v\otimes w-(-1)^{[v][w]}w\otimes v$ with homogeneous $v\in V$, $w\in W$. The isomorphism~\eqref{AotimesB} is induced by the embedding of $\NN_0$-graded super-algebras
\begin{align} \label{AotimesBEmb}
 TV\otimes TW=\bigoplus_{k,l\ge0}V^{\otimes k}\otimes W^{\otimes l}\hookrightarrow T(V\oplus W).
\end{align}
The isomorphism~\eqref{AcircB} is induced by the graded isomorphism
\begin{align}
 \phi\colon TV\circ TW=\bigoplus_{k\ge0}V^{\otimes k}\otimes W^{\otimes k}\isoright \bigoplus_{k\ge0}(V\otimes W)^{\otimes k}=T(V\otimes W), \label{phiT} \\
 \phi_k\colon v_1\otimes\cdots\otimes v_k\otimes w_1\otimes\cdots\otimes w_k\mapsto(-1)^{\bar p}(v_1\otimes w_1)\otimes\cdots\otimes(v_k\otimes w_k), \label{phik}
\end{align}
where $\bar p=\sum\limits_{1\le t<s\le k}[v_s][w_t]$.
\end{Prop}

\noindent{\bf Proof.} First note that any element of $T(V\oplus W)$ is equal to a sum of elements of the spaces $V^{\otimes k}\otimes W^{\otimes l}$ modulo the ideal generated by $[V,W]$, so the embedding~\eqref{AotimesBEmb} gives the isomorphism $TV\otimes TW\cong T(V\oplus W)/([V,W])$ in $\GrSVect{\NN_0}$. One can check that it respects the multiplication: e.g. $w\cdot v=(-1)^{[v][w]}v\otimes w\mapsto(-1)^{[v][w]}v\otimes w=w\otimes v$ $\mod[V,W]$. Due to Lemma~\ref{LemTPQ} the super-vector space $\A\otimes\B$ is isomorphic the quotient of $TV\otimes TW$ over $(R)\otimes TW+TV\otimes(S)$, so we obtain~\eqref{AotimesB}.
Further, it is straightforward to prove that the formula~\eqref{phik} defines a morphism of $\NN_0$-graded super-algebras~\eqref{phiT}. By virtue of Lemma~\ref{LemTPQ} the kernel of the epimorphism $TV\circ TW\twoheadrightarrow\A\circ\B$ is the ideal $(R)\circ TW+TV\circ(S)$. The isomorphism~\eqref{phiT} maps this ideal to $(R_{\mathbf w})$. \qed

A quadratic super-algebra $\A\cong TV/(R)$ is {\it purely even} or {\it purely odd} iff $V=V_\bz$ or $V=V_\bu$ respectively (see~\cite[\S~4.1.4]{Sqrt}). The purely even ones form the subcategory $\QA$ of the quadratic algebras $\A=\A_\bz$. The purely odd ones are such that $\A_\bz=\bigoplus\limits_{k\ge0}\A_{2k}$, $\A_\bu=\bigoplus\limits_{k\ge0}\A_{2k+1}$. Recall that Manin defined in~\cite{Manin87}, \cite{Manin88} four binary operations $\circ$, $\bullet$, $\otimes$, $\ootimes$. The operations~\eqref{AotimesB} and \eqref{AcircB} generalises the Manin's binary operations $\otimes$ and $\circ$: they respectively coincide for purely even algebras. The operation~\eqref{AotimesB} for purely odd algebras coincides with the operation $\ootimes$. Now let us generalise the {\it black Manin product} $\bullet$ and "odd" tensor product $\ootimes$ (from the purely even case). For $\A=TV/(R)$ and $\B=TW/(S)$ as in Proposition~\ref{Propotimescirc} define
\begin{align}
 &\A\bullet\B=T(V\otimes W)/(R_{\mathrm b}), &&R_{\mathrm b}=\sigma^{(23)}(R\otimes S); \label{AbulletB} \\
 &\A\ootimes\B=T(V\oplus W)/(R''), &&R''=R\oplus S\oplus[V,W]_+, \label{AootimesB}
\end{align}
where $[V,W]_+\subset V\otimes W\oplus W\otimes V$ is spanned by $v\otimes w+(-1)^{[v][w]}w\otimes v$.

\begin{Prop}
 The formula~\eqref{AbulletB} gives a symmetric monoidal category $(\QSA,\bullet)$ with unit object $\KK[\varepsilon]/(\varepsilon^2)\in\QA$ and symmetric structure~\eqref{sigmaSuper}. The functor $\bullet$ on morphisms is uniquely defined by the condition $(f\bullet g)_1=f_1\otimes g_1$.
\end{Prop}

\noindent{\bf Proof.} Let $f\colon\A\to\wt\A$ and  $g\colon\B\to\wt\B$, where $\wt\A=T\wt V/(\wt R)$, $\wt\B=T\wt W/(\wt S)$. Then we have $(f_1\otimes g_1\otimes f_1\otimes g_1)(R_{\mathbf b})=\sigma^{(23)}(f_1\otimes f_1\otimes g_1\otimes g_1)(R\otimes S)=\sigma^{(23)}(\wt R\otimes\wt S)$. \qed

There is a natural epimorphism $\A\bullet\B\twoheadrightarrow\A\circ\B$ whose first degree component is the map $\id\colon\A_1\otimes\B_1\to\A_1\otimes\B_1$.

The subcategory $\FQSA\subset\QSA$ is monoidal with respect to all three products $\otimes$, $\circ$, $\bullet$. 

\bb{Coproduct of quadratic super-algebras.}
Let us obtain the categorical coproduct of $\A=TV/(R)$ and $\B=TW/(S)$ in $\QSA$.

\begin{Prop}
 The coproduct of $\A$ and $\B$ is
\begin{align}
 \A\amalg\B=T(V\oplus W)/(R\oplus S).
\end{align}
It gives the symmetric monoidal categories $(\QSA,\amalg)$ and $(\FQSA,\amalg)$ with the unit object $\KK$.
\end{Prop}

\noindent{\bf Proof.} One can prove that $\A\amalg\B$ is the coproduct in $\QSA$ by a direct generalisation of~\cite[\S~4.1.9]{Sqrt} to the super-case. The initial object is $\KK$. Since $\KK\in\FQSA$ and a coproduct of $\A,\B\in\FQSA$ in $\QSA$ belongs to $\FQSA$, it coincides with the coproduct in $\FQSA$. As any category with finite coproducts the categories $\QSA$ and $\FQSA$ are symmetric monoidal with respect to $\amalg$. \qed

We have the natural epimorphism
\begin{math}
 \A\amalg\B\twoheadrightarrow\A\otimes\B.
\end{math} 

\subsection{Some functors}
\label{sec32}

\bb{Dualisation functor.} Let $V\in\SVect$. Its dual $V^*=\bfhom(V,\KK)$ is the super-vector space with the components $V^*_\bz=\{\lambda\colon V\to\KK\mid \lambda(V_\bu)=0\}$, $V^*_\bu=\{\lambda\colon V\to\KK\mid \lambda(V_\bz)=0\}$. The dualisation can be considered as the contravariant functor $(-)^*\colon\SVect\to\SVect$, $V\mapsto V^*$. It translates a morphism $f\colon V\to W$ to $f^*\colon W^*\to V^*$, where $f^*(\mu)(v)=\mu\big(f(v)\big)$, $\mu\in W^*$. The monoidal structure of this functor is defined by the embedding
\begin{align} \label{Dualphi}
 &V^*\otimes W^*\hookrightarrow(V\otimes W)^*, &&(\lambda\otimes\mu)(v\otimes w)=(-1)^{[\mu][v]}\lambda(v)\mu(w).
\end{align}
This gives a contravariant symmetric {\bf lax} monoidal functor $(-)^*\colon(\SVect,\otimes)\to(\SVect,\otimes)$. In the finite-dimensional case we obtain a contravariant symmetric {\bf strong} monoidal functor $(-)^*\colon(\FSVect,\otimes)\to(\FSVect,\otimes)$. More generally, the formula~\eqref{Dualphi} gives the isomorphism $V^*\otimes W^*\cong (V\otimes W)^*$ natural in $V\in\FVect$ and $W\in\Vect$.

The dualisation operation has also a structure of contravariant symmetric {\bf strong} monoidal functors $(-)^*\colon(\SVect,\oplus)\to(\SVect,\oplus)$ and $(-)^*\colon(\FSVect,\oplus)\to(\FSVect,\oplus)$.

\begin{Lem} \label{Lembot}
 Let $V_0$ be a super-subspace of a super-vector space $V\in\SVect$. Consider the projection $\pi\colon V\to V/V_0$. The image of $\pi^*\colon(V/V_0)^*\to V^*$ equals
\begin{align}
 V_0^\bot:=\{\lambda\in V^*\mid\lambda(v_0)=0\;\forall\,v_0\in V_0\}.
\end{align}
It is a super-subspace $V_0^\bot\subset V^*$ isomorphic to $(V/V_0)^*$ via $\pi^*$. Suppose $V,W\in\FSVect$ and let $W_0\subset W$ be also a super-subspace, then $V_0^\bot\otimes W_0^\bot$ coincide with $(V_0\otimes W+V\otimes W_0)^\bot$ as a super-subspace of $V^*\otimes W^*=(V\otimes W)^*$.
\end{Lem}

\noindent{\bf Proof.} It is obvious that $\Ker(\pi^*)=0$ and $\Im(\pi^*)\subset V_0^\bot$. Let $t\colon V_0^\bot\to(V/V_0)^*$ be the map $t(\lambda)(\bar v)=\lambda(v)$, where $\lambda\in V_0^\bot$, $v\in V$, $\bar v\in V/V_0$. Since $\pi^*\big(t(\lambda)\big)=\lambda$ $\;\forall\,\lambda\in V_0^\bot$, we obtain $\Im(\pi^*)=V_0^\bot$, so $\pi^*$ induces the isomorphism $(V/V_0)^*\cong V_0^\bot$. By applying the functor $(-)^*$ to~\eqref{VV0WW0} we obtain $V_0^\bot\otimes W_0^\bot=(V_0\otimes W+V\otimes W_0)^\bot$. \qed

\begin{Lem} \label{Lemxibot}
 Let $V_0$ and $W_0$ be super-subspaces of $V\in\FVect$ and $W\in\Vect$ respectively. Identify the operators $\xi\in\bfhom(V,W)$ with the elements $\xi\in V^*\otimes W$. Then the condition $\xi(V_0)\subset W_0$ is equivalent to $\xi\in V_0^\bot\otimes W+V^*\otimes W_0$.
\end{Lem}

\noindent{\bf Proof.} The first condition means exactly that $\xi$ belongs to the kernel of the epimorphism $p\colon V^*\otimes W\twoheadrightarrow V_0^*\otimes(W/W_0)$. Note that $V_0=(V_0^\bot)^\bot\cong(V^*/V_0^\bot)^*$ and hence $V_0^*\cong(V^*/V_0^\bot)$. By applying Lemma~\ref{LemTPQ} for $V_0^\bot\subset V^*$ and $W_0\subset W$ we see that $\Ker p=V_0^\bot\otimes W+V^*\otimes W_0$. \qed

\bb{Functors of parity change: $\Pi$ and $\wh\Pi$.} Denote by $\Pi\colon\SVect\to\SVect$ the functor changing parity: $(\Pi V)_\bz=V_\bu$, $(\Pi V)_\bu=V_\bz$. This is an autoequivalence of $\SVect$, it has symmetric strong monoidal structure $(\SVect,\oplus)\to(\SVect,\oplus)$. It preserves the tensor product, but not as a monoidal functor: the isomorphism $(\Pi V)\otimes(\Pi W)\cong V\otimes W$ is natural in $V,W\in\SVect$. It commutes with $(-)^*$ in the sense of natural isomorphism $(\Pi V)^*=\Pi(V^*)$. Note that the composition $\Pi V^*\otimes\Pi W^*\cong(\Pi V\otimes\Pi W)^*\cong(V\otimes W)^*\cong V^*\otimes W^*$ coincides with $\Pi V^*\otimes\Pi W^*\xisoright{\sigma}\Pi W^*\otimes\Pi V^*\cong W^*\otimes V^*\xisoright{\sigma}V^*\otimes W^*$, $\lambda\otimes\mu\mapsto-(-1)^{[\lambda]+[\mu]}\lambda\otimes\mu$, where $\lambda\in V^*$, $\mu\in W^*$.

Define the functor $\wh\Pi\colon\QSA\to\QSA$ by $\wh\Pi\big(TV/(R)\big)=T\Pi V/(R)$. This is an autoequivalence of $\QSA$, more precisely, $\wh\Pi^2=\id$. It maps the purely even quadratic super-algebras to the purely odd ones and vice versa. By restricting it to the subcategory $\QA\subset\QSA$ we obtain the functor $\wh\Pi\colon\QA\to\QSA$, which gives the equivalence of the subcategories of purely even and odd quadratic super-algebras.

\begin{Prop}
 The functor $\wh\Pi\colon\QSA\to\QSA$ has a structure of symmetric {\bf strong} monoidal functor $\wh\Pi\colon(\QSA,\amalg)\to(\QSA,\amalg)$. The restricted functor $\wh\Pi\colon\QA\to\QSA$ has symmetric {\bf strong} monoidal structures $\wh\Pi\colon(\QA,\otimes)\to(\QSA,\ootimes)$ and $\wh\Pi\colon(\QA,\ootimes)\to(\QSA,\otimes)$, i.e. $\wh\Pi\A\otimes\wh\Pi\B=\wh\Pi(\A\ootimes\B)$, $\wh\Pi\A\ootimes\wh\Pi\B=\wh\Pi(\A\otimes\B)$ for any $\A,\B\in\QA$. We also have $\wh\Pi\A\otimes\B=\wh\Pi(\A\otimes\wh\Pi\B)$, $\wh\Pi\A\ootimes\B=\wh\Pi(\A\ootimes\wh\Pi\B)$, $\wh\Pi\A\circ\wh\Pi\B=\A\circ\B$, $\wh\Pi\A\bullet\wh\Pi\B=\A\bullet\B$, $\wh\Pi\A\circ\B=\wh\Pi(\A\circ\B)$, $\wh\Pi\A\bullet\B=\wh\Pi(\A\bullet\B)$ $\;\forall\,\A,\B\in\QA$.
\end{Prop}

\noindent{\bf Proof.} Let $\A=TV/(R)$ and $\B=TW/(S)$ be quadratic super-algebras. We have the identifications $\wh\Pi\A\amalg\wh\Pi\B=T(\Pi V\oplus\Pi W)/(R\oplus W)=\wh\Pi(\A\amalg\B)$. If $\A,\B\in\QA$, then $V,W\in\Vect$ and hence $[\Pi V,\Pi W]=[V,W]_+$. This implies $\wh\Pi\A\otimes\wh\Pi\B=T(\Pi V\oplus\Pi W)/(R\oplus S\oplus[\Pi V,\Pi W])=\wh\Pi(\A\ootimes\B)$. The next three isomorphisms follows from $[\Pi V,\Pi W]_+=[V,W]$, $[\Pi V,W]=\Pi[V,W]=[V,\Pi W]$, $[\Pi V,W]_+=\Pi[V,W]_+=[V,\Pi W]_+$. Further, we derive $\wh\Pi\A\bullet\wh\Pi\B=T(\Pi V\otimes\Pi W)/\big(\sigma_{\Pi V,\Pi W}^{(23)}(R\otimes S)\big)=T(V\otimes W)/\big(\sigma_{V,W}^{(23)}(R\otimes S)\big)=\A\bullet\B$, since $\sigma_{\Pi V,\Pi W}^{(23)}(r\otimes s)=-\sigma_{V,W}^{(23)}(r\otimes s)$ for any $r\in V^{\otimes2}$ and $s\in W^{\otimes2}$. The other isomorphisms are checked similarly. \qed


\bb{Functors $(-)_1$ and $(-)_1^*$.} For each $k\in\ZZ$ we have the functor $(-)_k\colon\GrSAlg{\ZZ}\to\SVect$ which extract the $k$-th component $\A_k\in\SVect$ from a $\ZZ$-graded super-algebra $\A$. Its value on the graded homomorphism $f\colon\A\to\B$ is the component $f_k\colon\A_k\to\B_k$. The most interesting functor is $(-)_1\colon\GrSAlg{\ZZ}\to\SVect$. Its restrictions $(-)_1\colon\QSA\to\SVect$ and $(-)_1\colon\FQSA\to\FSVect$ have structures of symmetric {\bf strong} monoidal functors
\begin{align}
 &(-)_1\colon(\QSA,\otimes)\to(\SVect,\oplus), &
 &(-)_1\colon(\FQSA,\otimes)\to(\FSVect,\oplus), \\
 &(-)_1\colon(\QSA,\amalg)\to(\SVect,\oplus), &
 &(-)_1\colon(\FQSA,\amalg)\to(\FSVect,\oplus), \\
 &(-)_1\colon(\QSA,\circ)\to(\SVect,\otimes), &
 &(-)_1\colon(\FQSA,\circ)\to(\FSVect,\otimes), \label{circ1} \\
 &(-)_1\colon(\QSA,\bullet)\to(\SVect,\otimes), &
 &(-)_1\colon(\FQSA,\bullet)\to(\FSVect,\otimes). \label{bullet1}
\end{align}

The composition of $(-)_1\colon\QSA\to\SVect$ with the dualisation $(-)^*\colon\SVect\to\SVect$ gives the contravariant functor $(-)_1^*\colon\QSA\to\SVect$, $\A\mapsto\A_1^*$. These functors commutes with the change of parity as follows: $(\wh\Pi\A)_1=\Pi\A_1$, $(\wh\Pi\A)_1^*=\Pi\A_1^*$, these isomorphisms are natural in $\A\in\QSA$.

The restricted functor $(-)_1^*\colon\FQSA\to\FSVect$ has structures of contravariant symmetric {\bf strong} monoidal functors
\begin{align}
 &(-)_1^*\colon(\FQSA,\otimes)\to(\FSVect,\oplus), &
 &(-)_1^*\colon(\FQSA,\amalg)\to(\FSVect,\oplus), \\
 &(-)_1^*\colon(\FQSA,\circ)\to(\FSVect,\otimes), &
 &(-)_1^*\colon(\FQSA,\bullet)\to(\FSVect,\otimes).
\end{align}

\bb{Functors $T$ and $T^*$.} Let us consider the functor $T\colon\SVect\to\QSA$ making the tensor algebra $TV\in\QSA$ for a super-vector space $V$. For a $\ZZ_2$-graded linear map $f\colon V\to W$ the there exists a unique morphism $Tf\colon TV\to TW$ with first degree component $(Tf)_1=f$. This means exactly that $T\colon\SVect\to\QSA$ is a left adjoint to $(-)_1\colon\QSA\to\SVect$. We have isomorphism $T(\Pi V)=\wh\Pi(TV)$ natural in $V\in\SVect$. The functor $T\colon\SVect\to\QSA$ is symmetric {\bf strong} monoidal as $(\SVect,\oplus)\to(\QSA,\amalg)$ and $(\SVect,\otimes)\to(\QSA,\circ)$ and symmetric {\it colax} monoidal as $(\SVect,\oplus)\to(\QSA,\otimes)$.

Denote the composition of the functors $(-)^*\colon\SVect\to\SVect$ and  $T\colon\SVect\to\QSA$ by $T^*\colon\SVect\to\QSA$. This functor as well as its restriction $T^*\colon\FSVect\to\FQSA$ is a faithful contravariant functor. Note that a composition of a (symmetric) strong monoidal functor with a (symmetric) lax/colax/strong monoidal functor is a (symmetric) lax/colax/strong monoidal functor. Hence we obtain the following contravariant symmetric monoidal functors:
\begin{align}
&\text{{\bf colax} monoidal} &&T^*\colon(\FSVect,\oplus)\to(\FQSA,\otimes), \label{laxTstFVectplus} \\
&\text{{\bf strong} monoidal} &&T^*\colon(\FSVect,\oplus)\to(\FQSA,\amalg), \label{TstFVectplus} \\
&\text{{\bf strong} monoidal} &&T^*\colon(\FSVect,\otimes)\to(\FQSA,\circ). \label{TstFVectotimes} 
\end{align}

\bb{Functors $S$ and $S^*$.} Recall two functors $S\colon\Vect\to\QA$ and $\Lambda\colon\Vect\to\QA$ which makes the symmetric algebra $SV$ and external algebra $\Lambda V$ from a vector space $V$. The functor $S$ can be extended to the super-case as follows: for $V\in\SVect$ we define $SV=SV_\bz\otimes\wh\Pi(\Lambda V_\bu)$, where $SV_\bz$ and $\Lambda V_\bu$ are considered as purely even quadratic super-algebras. Note that $SV$ is a quadratic super-algebra with the relations $xy=(-1)^{[x][y]}yx$, where $x,y\in V$ are homogeneous. We obtain the functor $S\colon\SVect\to\QSA$, which translates a $\ZZ_2$-graded map $f$ to $Sf_\bz\otimes\wh\Pi(\Lambda f_\bu)$.

Consider the full subcategory $\CommQSA\subset\QSA$ consisting of $\A\in\QSA$ which are commutative as super-algebras: $ab=(-1)^{[a][b]}ba$. Note that $SV\in\CommQSA$ for any $V\in\SVect$. Hence $S$ can be considered as a functor $\SVect\to\CommQSA$, which is a left adjoint to $(-)_1\colon\CommQSA\to\SVect$. Alternatively, one can say that $(-)_1\colon\CommQSA\to\SVect$ is a right adjoint for $S\colon\SVect\to\QSA$ relative to the full subcategory $\CommQSA\subset\QSA$.

By composing $(-)^*\colon\SVect\to\SVect$ and $S\colon\SVect\to\QSA$ we obtain a contravariant faithful functors $S^*\colon\SVect\to\QSA$ and $S^*\colon\FSVect\to\FQSA$.

\begin{Prop}
 The functors
 $S\colon(\SVect,\oplus)\to(\QSA,\otimes)$, $S\colon(\FSVect,\oplus)\to(\FQSA,\otimes)$
are symmetric {\bf strong} monoidal. Hence we have contravariant
\begin{align}
 &\text{symmetric {\bf strong} monoidal functor} &S^*\colon(\FSVect,\oplus)\to(\FQSA,\otimes).
\end{align}
The functors
 $S\colon(\SVect,\otimes)\to(\QSA,\circ)$, $S\colon(\FSVect,\otimes)\to(\FQSA,\circ)$
are symmetric {\bf colax} monoidal. Hence we have contravariant
\begin{align}
 &\text{ symmetric {\bf colax} monoidal functor} &&S^*\colon(\FSVect,\otimes)\to(\FQSA,\circ).
\end{align}
\end{Prop}

\noindent{\bf Proof.} Proposition~\ref{Propotimescirc} implies that the morphisms~\eqref{AotimesBEmb} and \eqref{phiT} induce the isomorphism $S(V\otimes W)\isoright SV\otimes SW$ and epimorphism $S(V\otimes W)\twoheadrightarrow SV\circ SW$ for $V,W\in\SVect$. \qed

\bb{Koszul duality $(-)^!$ on $\FQSA$.} Let us consider a super-version of the Manin's contravariant functor $(-)^!\colon\FQA\to\FQA$ introduced in~\cite{Manin87}. For a quadratic super-algebra $\A=TV/(R)\in\FQSA$ we set
\begin{align} \label{DefKD}
 &\A^!=TV^*/(R^\bot),  &&R^\bot=\{\xi\in V^*\otimes V^*\mid\xi(r)=0\;\forall\,r\in R\},
\end{align}
where $\xi(r)$ is defined by the formula~\eqref{Dualphi}. Due to Lemma~\ref{Lembot} the set $R^\bot$ is a super-subspace of $V^*\otimes V^*=(V\otimes V)^*$, so $\A^!\in\FQSA$. For any morphism $f\colon TV/(R)\to TW/(S)$ in $\FQSA$ we have $(f_1^*\otimes f_1^*)S^\bot\subset R^\bot$, hence we obtain the contravariant functor $(-)^!\colon\FQSA\to\FQSA$ such that
\begin{align} \label{A1KD}
 &(\A^!)_1=\A_1^*, &&(f^!)_1=f_1^*, &&\text{and} &&(\A^!)^!=\A, &&(f^!)^!=f.
\end{align}
It commutes with $\wh\Pi\colon\FQSA\to\FQSA$ in the sense that we have a natural isomorphism $(\wh\Pi\A)^!\cong\wh\Pi(\A^!)$. The right property~\eqref{A1KD} implies that the functor $(-)^!\colon\FQSA\to\FQSA$ is fully faithful.

\begin{Prop}
 The contravariant functor $(-)^!$ has symmetric strong monoidal structures
\begin{align*}
 &(\FQSA,\circ)\to(\FQSA,\bullet), &
 &(\FQSA,\bullet)\to(\FQSA,\circ), \\
 &(\FQSA,\otimes)\to(\FQSA,\ootimes) &
 &(\FQSA,\ootimes)\to(\FQSA,\otimes)
\end{align*}
given by the identifications $\KK[u]^!=\KK[\varepsilon]/(\varepsilon^2)$, $\big(\KK[\varepsilon]/(\varepsilon^2)\big)^!=\KK[u]$, $\KK^!=\KK$ and the isomorphisms
\begin{align}
 &(\A\circ\B)^!\cong\A^!\bullet\B^!, &
 &(\A\bullet\B)^!\cong\A^!\circ\B^!, \label{IsoKD1} \\
 &(\A\otimes\B)^!\cong\A^!\ootimes\B^!, &
 &(\A\ootimes\B)^!\cong\A^!\otimes\B^! \label{IsoKD2}
\end{align}
natural in $\A,\B\in\FQSA$.
\end{Prop}

\noindent{\bf Proof.} Let $\A=TV/(R),\B=TW/(S)\in\FQSA$. By means of Lemma~\ref{Lembot} we derive $R^\bot\otimes S^\bot=(R\otimes W^{\otimes2}+V^{\otimes2}\otimes S)^\bot=\sigma^{(23)}R_{\mathbf w}^\bot$, where $R_{\mathbf w}$ is defined by~\eqref{AcircB}. Then $\A^!\bullet\B^!=T(V^*\otimes W^*)/\big(\sigma^{(23)}(R^\bot\otimes S^\bot)\big)=T(V^*\otimes W^*)/(R_{\mathbf w}^\bot)\cong(\A\circ\B)^!$. The naturality is deduced by application of the faithful functor $(-)_1$ and using the properties~\eqref{A1KD}, $(\A\circ\B)_1=(\A\bullet\B)_1=\A_1\otimes\B_1$ (see~\eqref{circ1}, \eqref{bullet1}). The third isomorphism is obtained from the formulae $(R\oplus S)^\bot=R^\bot\oplus S^\bot$, $[V,W]^\bot=[V^*,W^*]_+$ and $(\A\otimes\B)_1=\A_1\oplus\B_1$. By substituting $\A\to\A^!$ and $\B\to\B^!$ and taking into account $(\A^!)^!=\A$ we obtain the right isomorphisms~\eqref{IsoKD1}, \eqref{IsoKD2}. \qed

\bb{Functors $\Lambda$ and $\Lambda^*$.} For a super-vector space $V\in\SVect$ define the quadratic algebra $\Lambda V=TV/(R)$, where $R\subset V^{\otimes2}$ is spanned by $x\otimes y+(-1)^{[x][y]}y\otimes x$ with homogeneous $x,y\in V$. We can decompose it to the purely even and purely odd parts as $\Lambda V=\Lambda V_\bz\ootimes\wh\Pi(SV_\bu)$. We obtain a functor $\Lambda\colon\SVect\to\QSA$, whose composition with $(-)_1$ is the identical functor. Define $\Lambda^*\colon\FSVect\to\FQSA$ as the composition of $(-)^*$ and $\Lambda\colon\FSVect\to\FQSA$.

\begin{Prop}
 The covariant functors
 $\Lambda\colon(\SVect,\oplus)\to(\QSA,\ootimes)$, $\Lambda\colon(\FSVect,\oplus)\to(\FQSA,\ootimes)$ and contravariant functor $\Lambda^*\colon(\FSVect,\oplus)\to(\FQSA,\ootimes)$
are fully faithful symmetric {\bf strong} monoidal.
The functors $\Lambda\colon(\SVect,\otimes)\to(\QSA,\bullet)$, $\Lambda\colon(\FSVect,\otimes)\to(\FQSA,\bullet)$, $\Lambda^*\colon(\FSVect,\otimes)\to(\FQSA,\bullet)$ are symmetric {\bf lax} monoidal. There are natural isomorphisms $(\Lambda V)^!\cong S^*(V)$, $(SV)^!\cong\Lambda^*(V)$.
\end{Prop}

\noindent{\bf Proof.} Straightforward. \qed

\bb{Functor $(-)^\op$.} \label{bbFunctop}
For a super-algebra $\A$ with the multiplication map $\mu\colon\A\otimes\A\to\A$ denote by $\A^\op$ the algebra with the opposite multiplication: $\mu^\op=\mu\cdot\sigma\colon\A\otimes\A\to\A$, $a\otimes b\mapsto(-1)^{[a][b]}ba$. If $\A$ is $\AA$-graded, then the opposite super-algebra $\A^\op$ is also $\AA$-graded. For a quadratic algebra $\A\cong TV/(R)$ we have $\A^\op\cong TV/(\sigma R)$, so we get an autoequivalence $(-)^\op\colon\QSA\to\QSA$. It commutes with the functor $\wh\Pi\colon\QSA\to\QSA$ and preserves the products $\otimes$, $\ootimes$, $\circ$, $\bullet$ and $\amalg$ as symmetric strong monoidal functor. The restricted functor $(-)^\op\colon\FQSA\to\FQSA$ commutes with $(-)^!\colon\FQSA\to\FQSA$.

\subsection{Internal cohom-functor}
\label{sec33}

\bb{Internal hom and cohom in $\FQSA$.} For finitely generated quadratic algebras the internal hom and cohom for the products $\bullet$ and $\circ$ respectively were obtained by Manin in~\cite{Manin87}, \cite{Manin88}. He derived a bijection
\begin{align} \label{HomABC}
 \Hom(\A\bullet\B,\C)\cong\Hom(\A,\B^!\circ\C)
\end{align}
natural in $\A,\B,\C\in\FQA$. This implies that the symmetric monoidal categories $(\FQA,\bullet)$ and $(\FQA,\circ)$ are closed and coclosed respectively. Let us generalise this fact for the quadratic super-algebras. First we consider two particular cases $\A=\KK[\varepsilon]/(\varepsilon^2)$ and $\C=\KK[u]$.

\begin{Lem} \label{LemHomABK}
 For $V,W\in\SVect$ we have the identification $\Hom(\KK,V^*\otimes W)=\Hom(V,W)$ and $\Hom(W\otimes V^*,\KK)=\Hom(W,V)$. It induces the bijections
\begin{align} \label{HomABK}
 &\Hom\big(\KK[\varepsilon]/(\varepsilon^2),\A^!\circ\B\big)\cong\Hom(\A,\B), &
 &\Hom\big(\B\bullet\A^!,\KK[u]\big)\cong\Hom(\B,\A)
\end{align}
natural in $\A\in\FQSA$, $\B\in\QSA$.
\end{Lem}

\noindent{\bf Proof.} By taking the even components of the natural isomorphisms of super-vector spaces $\bfhom(\KK,V^*\otimes W)=V^*\otimes W=\bfhom(V,W)$, $\bfhom(W\otimes V^*,\KK)=W^*\otimes V=\bfhom(W,V)$ and forgetful functor $\Vect\to\Set$ we obtain the identifications of the corresponding Hom-sets.
 Let $\A=TV/(R)$, $\B=TW/(S)$, where $V\in\FVect$, $W\in\Vect$. The elements of the Hom-sets in the left formula~\eqref{HomABK} are morphisms in $\QSA$ uniquely defined by their first degree components $\lambda\colon\KK\to V^*\otimes W$ and $f\colon V\to W$ which satisfy the conditions $(\lambda\otimes\lambda)\in\sigma^{(23)}(R^\bot\otimes W^{\otimes2}+(V^*)^{\otimes2}\otimes S)$ and $(f\otimes f)R\subset S$.
If we identify $\lambda$ with $f$ via $\Hom(\KK,V^*\otimes W)=\Hom(V,W)$, then we have $\lambda=\sum_i \lambda_i\otimes w_i$, $f=\sum_i w_i\lambda_i$ for some $\lambda_i\in V^*$, $w_i\in W$ such that $[\lambda_i]=[w_i]$. Note that the map $f\otimes f\colon V\otimes V\to W\otimes W$ acts on $v\otimes v'\in V\otimes V$ as $(f\otimes f)(v\otimes v')=f(v)f(v')=\sum_{i,j}(w_i\otimes w_j)\lambda_i(v)\lambda_j(v')$, hence $f\otimes f=\sum_{i,j}(-1)^{[\lambda_i][\lambda_j]}(w_i\otimes w_j)(\lambda_i\otimes\lambda_j)$. By applying Lemma~\ref{Lemxibot} to the even operator $\xi=\sigma^{(23)}(\lambda\otimes\lambda)=\sum_{i,i}(-1)^{[\lambda_i][\lambda_j]}\lambda_i\otimes\lambda_j\otimes w_i\otimes w_j\in \big((V^{\otimes2})^*\otimes W^{\otimes2}\big)_\bz$ we see that these two conditions are equivalent to each other. The right isomorphism~\eqref{HomABK} is proved similarly (for $\B\in\FQSA$ it follow from the left isomorphism~\eqref{HomABK}, the isomorphism~\eqref{IsoKD1} and the fact that $(-)^!\colon\FQSA\to\FQSA$ is fully faithful). \qed

\begin{Prop}
 The identification $\Hom(\A_1\otimes\B_1,\C_1)=\Hom(\A,\B_1^*\otimes\C)$ induces the bijection $\Hom(\A\bullet\B^!,\C)\cong\Hom(\A,\B\circ\C)$ (or, equivalently,~\eqref{HomABC})
natural in $\A,\B\in\FQSA$, $\C\in\QSA$. It can be also regarded as a bijection natural in $\A\in\QSA$, $\B,\C\in\FQSA$. As consequence, the symmetric monoidal categories $(\QSA,\bullet)$ and $(\QSA,\circ)$ are closed and coclosed respectively relative to the subcategory $\bfP=\FQSA$. The hom- and cohom functors have the form $\bfhom_{(\FQSA,\bullet)}(\B,\C)=\B^!\circ\C$ and $\bfcohom_{(\FQSA,\circ)}(\B,\A)=\A\bullet\B^!$.
\end{Prop}

\noindent{\bf Proof.} Lemma~\ref{LemHomABK} and the isomorphisms~\eqref{IsoKD1} give the bijections
 $\Hom(\A\bullet\B^!,\C)\cong
\Hom\big(\KK[\varepsilon]/(\varepsilon^2),(\A\bullet\B^!)^!\circ\C\big)\cong\Hom\big(\KK[\varepsilon]/(\varepsilon^2),\A^!\circ\B\circ\C\big)\cong\Hom(\A,\B\circ\C)$ natural in $\A,\B\in\FQSA$, $\C\in\QSA$ and the bijections 
 $\Hom(\A\bullet\B,\C)\cong\Hom(\A\bullet\B\bullet\C^!,\KK[u])\cong\Hom(\A\bullet(\B^!\circ\C)^!,\KK[u])\cong\Hom(\A,\B^!\circ\C)$ natural in $\A\in\QSA$, $\B,\C\in\FQSA$.
 The hom- and cohom-functors are obtained from these bijections. \qed

\bb{Internal cohom for $(\FQSA,\circ)$ via the internal hom for $\SVect$.}
Let us derive more explicit expression for the cohom-functor generalising~\cite[Prop.~4.5]{Sqrt}. The first component is $\bfcohom(\B,\A)_1=W^*\otimes V=\bfhom(W,V)$, where $V=\A_1$, $W=\B_1$. There is a natural isomorphism
\begin{align} \label{bfhomWV2}
 \bfhom(W,V)^{\otimes2}=W^*\otimes V\otimes W^*\otimes V\xisoright{\sigma^{(23)}}W^*\otimes W^*\otimes V\otimes V=\bfhom(W^{\otimes2},V^{\otimes2}).
\end{align}

\begin{Prop} \label{Propcohom}
The cohom-functor for $\A=TV/(R)$ and $\B=TW/(S)$ has the form
\begin{align}
 &\bfcohom(\B,\A)=\A\bullet\B^!=T\wt V/(\wt R), &
 &\text{where\quad}\wt V=\bfhom(W,V)
\end{align}
and $\wt R\subset\wt V^{\otimes2}$ is the preimage of $\big\{\xi\in\bfhom(W^{\otimes2},V^{\otimes2})\mid\xi(W^{\otimes 2})\subset R, \xi(S)=0\big\}$ under~\eqref{bfhomWV2}. For two morphisms $f$ in $\FQSA$ and $g$ in $\QSA$ the morphism $\bfcohom(f,g)$ is defined by the component $\bfcohom(f,g)_1=\bfhom(f_1,g_1)=f_1^*\cdot(g_1)_*$.
\end{Prop}

\noindent{\bf Proof.} By the definitions~\eqref{AbulletB} and~\eqref{DefKD} we have $\A\bullet\B^!=T\wt V/(\wt R)$ for the super-vector space $\wt V=V\otimes W^*=\bfhom(W,V)$ and super-subspace $\wt R=\sigma^{(23)}(S^\bot\otimes R)\subset\wt V^{\otimes2}$. If $\xi$ belongs to $\sigma^{(23)}\wt R=S^\bot\otimes R$, then $\xi(W^{\otimes2})\subset R$ and $\xi(S)=0$. Conversely, due to Lemma~\ref{Lemxibot} the latter two conditions imply $\xi\in(W^{\otimes2})^*\otimes R$ and $\xi\in S^\bot\otimes V^{\otimes2}$ respectively, hence $\xi$ belongs to the intersection $\big((W^{\otimes2})^*\otimes R\big)\cap\big(S^\bot\otimes V^{\otimes2}\big)=S^\bot\otimes R$. \qed

Note that the map $\sigma\colon\A_1\otimes\B_1^*\isoright\B_1^*\otimes\A_1$ gives the natural isomorphism
\begin{align} \label{BAAB}
 \bfcohom(\B,\A)=\A\bullet\B^!\cong\B^!\bullet\A=\bfcohom(\A^!,\B^!).
\end{align}

According to~\cite[Th.~4.3]{Sgrt} the monoidal categories $(\GrAlg{\NN_0},\circ)$ and $(\GrAlg{\ZZ},\circ)$ are coclosed relative to $\FQA$. This result can be generalised to the super-case. We prove this below by using Manin matrices (see p.~\ref{bbUMM}).

\section{Manin matrices and quantum representations}
\label{sec4}

Now we generalise results on Manin matrices~\cite{Smm} and on quantum representations~\cite{Sqrt} to the super-case.

\subsection{Manin matrices for quadratic super-algebras}
\label{sec41}

\bb{Format.}
Consider a basis $e\colon I\to W$ of a super-vector space $W\in\SVect$. We always suppose that a basis is homogeneous: $e_i\in W_{\bar k_i}$, where $e_i=e(i)$, $i\in I$. The map $\bfk\colon I\to\ZZ_2$, $i\mapsto\bar k_i$, is called {\it format}%
\footnote{In~\cite{ManinG}, \cite{ManinBook91} Manin used the term "format" for a format of a matrix.
}
of the basis $(e_i)$. In other words, a format $\bfk$ is given by the decomposition $I=I_\bz\amalg I_\bu$ such that $[e_i]=\bz$ for $i\in I_\bz$ and $[e_i]=\bu$ for $i\in I_\bu$.

In the case $W\in\FSVect$ one can take $I=\{1,2,\ldots,d\}$. Then the format is a finite sequence $\bfk=(\bar k_1,\bar k_2,\ldots,\bar k_d)$. We can reorder a basis $(e_i)$ to the {\it standard} format $\bfk=(\bz,\ldots,\bz,\bu,\ldots,\bu)$, where the numbers of $\bz$ and $\bu$ are equal to $m=\dim W_\bz$ and $n=\dim W_\bu=d-m$.

Let $\bfk\colon I\to\ZZ_2$ and $\bfl\colon\wt I\to\ZZ_2$ be maps from the sets $I$ and $\wt I$. Consider a map $I\times\wt I\to R$, $(i,j)\mapsto M_{ij}=M^i_j\in R$. It is called {\it matrix $(M^i_j)$ of a format $\bfk\times\bfl$ over a super-vector space} $R\in\SVect$ if the {\it entries} $M_{ij}$ are homogeneous elements of degree $\bar k_i+\bar l_j$. For ordered $I$ and $\wt I$ (such as $I=\{1,2,\ldots,d\}$, $\wt I=\{1,2,\ldots,\wt d\}$) the matrix is written as a table
 $\begin{pmatrix}
 M^1_1& \ldots & M^1_{\wt d} \\
 \vdots &\ddots &\vdots \\
 M^d_1& \ldots & M^d_{\wt d}
\end{pmatrix}$. If the formats $\bfk$ and $\bfl$ are standard, then it has the block form
 $\begin{pmatrix}
 A & B \\
 C & D
\end{pmatrix}$,
where $A\in\Mat_{\wt m\times m}(R_\bz)$, $B\in\Mat_{\wt n\times m}(R_\bu)$, $C\in\Mat_{\wt m\times n}(R_\bu)$, $D\in\Mat_{\wt n\times n}(R_\bz)$.

\bb{Matrix of an operator.} \label{bbMatrOp}
Let $W\in\SVect$ and $\R\in\SAlg$. The tensor product $W\otimes\R$ is a free right $\R$-module. We identify $w=w\otimes1_\R$, so the element $w\otimes r\in\R\otimes W$ can be written as $wr$. The isomorphism $\sigma\colon W\otimes\R\isoright\R\otimes W$ gives the structure of the left module: $r(wr')=(-1)^{[w][r]}wrr'$, $r,r\in\R$, $w\in W$.

Let $\wt W\in\SVect$. The $\ZZ_2$-graded linear maps $M\colon\wt W\to W\otimes\R$ are in one-to-one correspondence with the $\ZZ_2$-graded morphisms of the right $\R$-modules $\alpha\colon\wt W\otimes\R\to W\otimes\R$ as $\alpha_M(\wt w\otimes r)=M(\wt w)r$, where $r\in\R$, $\wt w\in\wt W$. Similarly, such $M$ uniquely defines a morphism of the left $\R$-modules $\alpha^M\colon\R\otimes\wt W\to\R\otimes W$. We have $\alpha_M=\alpha^M$ iff $M(\wt w)r=rM(\wt w)$ $\;\forall\,\wt w\in\wt W, r\in\R$.

Any morphism $M\in\Hom(\wt W,W\otimes\R)$ can be identified with an even element of the left free $\R$-module $\R\otimes\bfhom(\wt W,W)$, i.e. with an operator $M\in\big(\R\otimes\bfhom(\wt W,W)\big)_\bz$.
Let $(e_i)_{i\in I}$ and $(\wt e_a)_{a\in\wt I}$ be bases of $W$ and $\wt W$ with formats $\bfk\colon I\to\ZZ_2$ and $\bfl\colon\wt I\to\ZZ_2$.
By {\it matrix of the operator} $M$ we mean the matrix $(M^a_i)$ defined by the formula $M\wt e_a=\sum_{i\in I}\wt e_iM^i_a $. This is a matrix of the format $\bfk\times\bfl$ over $\R$. For a fixed bases $(e_i)$ and $(\wt e_a)$ we identify the operator with its matrix: $M=(M^i_a)$. The operator $W\to W\otimes\R$ with entries $\delta^i_j$ we denote by $1$, it corresponds to the identity $\R$-module morphism $\id\colon W\otimes\R\to W\otimes\R$.

Let $W'\in\SVect$ has a basis $(e'_j)_{j\in I'}$ of a format $\bfk'\colon I'\to\ZZ_2$. Let $(N^a_j)$ be the matrix of an operator $N\colon W'\to\wt W\otimes\R$, $Ne'_j=\sum_{a\in\wt I}\wt e_aN^a_j$, it has the format $\bfl\times\bfk'$. Define the composition $MN$ as the operator $W'\to W\otimes\R$ corresponding to the composition $W'\otimes\R\xrightarrow{\alpha_{M}}\wt W\otimes\R\xrightarrow{\alpha_N} W\otimes\R$, i.e.
\begin{align}
 MN\colon W'\xrightarrow{N}\wt W\otimes\R\xrightarrow{M\otimes\id_\R}W\otimes\R\otimes\R\xrightarrow{\id\otimes\id\otimes\mu_\R}W\otimes\R,
\end{align}
where $\mu_\R\colon\R\otimes\R\to\R$, $\mu_\R(r\otimes r')=rr'$, is the multiplication in $\R$. The matrix of $MN$ has the format $\bfk\times\bfk'$, its entries are $(MN)^i_j=\sum_{a\in\wt I}M^i_aN^a_j$.

If $M$ and $N$ are operators over some super-subalgebras of $\R$, then we can regard them as operators over $\R$, so the composition $MN$ is defined as an operator over $\R$ in the same way.

For a right $\R$-module $R$ denote $R^*=\bfhom(R,\R)$. This is a contravariant functor from the category of right $\R$-modules to the category of left $\R$-modules. It translates free right modules to free left modules: $(W\otimes\R)^*=\R\otimes W^*$. The morphism $\alpha_M\colon\wt W\otimes\R\to W\otimes\R$ is translated to the morphism $\alpha^*_M\colon\R\otimes W^*\to\R\otimes\wt W^*$, $\alpha_M^*(r\mu)(\wt w r')=r\mu\big(\alpha_M(\wt w r')\big)=r\mu(M\wt w)r'$, where $r,r'\in\R$, $\mu\in W^*$, $\wt w\in\wt W$ and we identified $\mu=1_\R\otimes\mu$.  We have $\alpha^*_M=\alpha^{M^\st}$ for the operator $M^\st\in\Hom(W^*,\wt W^*\otimes\R)$ defined as $M^\st(\mu)(\wt w)=\mu(M\wt w)$. Let $(e^i)$ and $(\wt e^a)$ be the bases of $W^*$ and $\wt W^*$ dual to $(e_i)$ and $(\wt e_a)$, they have the same formats $\bfk$ and $\bfl$. One can check that $M^\st(e^i)=\sum_{a\in\wt I}M^i_a\wt e^a=\sum_{a\in\wt I}\wt e^a M^\st_{ai}$, where $(M^\st)^a_i=(-1)^{(\bar k_i+\bar l_a)\bar l_a}M^i_a$ are entries of the {\it super-transposed} matrix. Since the operation $(-)^\st$ is not involutive we need also the inverse of super-transposition. Define a matrix $M^\ist$ by the formulae $(M^\ist)^\st=M=(M^\st)^\ist$; explicitly, $(M^\ist)^a_i=(-1)^{(\bar k_i+\bar l_a)\bar k_i}M^i_a$.

Note that in the case $\R=\KK$ we have the operator $M\colon\wt W\to W$. The operator $M^\st$ coincides with the image of $M$ under the functor $(-)^*\colon\SVect\to\SVect$, i.e. $M^\st=M^*\colon W^*\to\wt W^*$.

An element $y\in\wt W\otimes\R$ has the unique decomposition $y=\sum_{a\in\wt I}\wt e_ay^a$, where $y^a\in\R$ are {\it right coordinates} of $y$. The morphism $\alpha_M\colon\wt W\otimes\R\to W\otimes\R$ has the form $\alpha_M(y)=\sum_{a\in\wt I}M(\wt e_a)y^a=\sum_{i\in I}e_ix^i$, where $x^i=\sum_{a\in\wt I}M^i_a y^a$ are right coordinates of $x=\alpha_M(y)$. If we consider $x$ and $y$ as column-vectors with entries $x^i$ and $y^a$, then $\alpha_M$ is equivalent to the action of the matrix $M$ from the left: $x=My$. Similarly the map $\alpha^*_M\colon\R\otimes W^*\to\R\otimes\wt W^*$ is equivalent to the action of the matrix $M$ from the right: $\phi=\psi M$, where $\psi=\sum_{i\in I}\psi_i e^i\in\R\otimes W^*$, $\phi=\alpha^*_M(\psi)=\sum_{a\in\wt I}\phi_a\wt e^a\in\R\otimes\wt W^*$ are considered as row-vectors with entries $\psi_i\in\R$ and $\phi_a=\sum_{i\in I}\psi_i M^i_a\in\R$ ({\it left coordinates} of $\psi$ and $\phi$).

Consider the tensor product $W\otimes\wt W$. Its basis $(e_i\otimes\wt e_a)$ has the format $\bfk\oplus\bfl\colon\{(i,a)\mid i\in I,a\in\wt I\}\to\ZZ_2$, $(i,a)\mapsto\bar k_i+\bar l_a$. Introduce the operators $M^{(1)}$ and $M^{(2)}$ as
\begin{align}
 &M^{(1)}\colon \wt W\otimes W'\xrightarrow{M\otimes\id_{W'}} W\otimes\R\otimes W'\xisoright{\sigma^{(23)}} W\otimes W'\otimes\R, \label{M1} \\
 &M^{(2)}\colon W'\otimes\wt W\xrightarrow{\id_{W'}\otimes M} W'\otimes W\otimes\R \label{M2}
\end{align}
(these notations depends on $W'$).
They have the entries
\begin{align}
 &(M^{(1)})^{ij}_{al}=(-1)^{\bar k'_j(\bar k_i+\bar l_a)}M^i_a\delta^j_l, &&(M^{(2)})^{ji}_{la}=\delta^j_lM^i_a.
\end{align}
Let $N\colon \wt W'\to\R\otimes W'$. The composition of the operators $M^{(1)}\colon\wt W\otimes\wt W'\to W\otimes\wt W'\otimes\R$ and $N^{(2)}\colon \wt W\otimes\wt W'\to\wt W\otimes W'\otimes\R$ is the operator
\begin{multline}
 M^{(1)}N^{(2)}\colon W\otimes\wt W'\xrightarrow{M\otimes N}W\otimes\R\otimes\wt W'\otimes\R\xrightarrow{\sigma^{(23)}}\\
 W\otimes W'\otimes\R\otimes\R\xrightarrow{\id\otimes\id\otimes\mu_\R}W\otimes W'\otimes\R. \label{M1N2}
\end{multline}
Its entries are $(M^{(1)}N^{(2)})^{ij}_{ab}=(-1)^{(\bar k_i+\bar l_a)\bar k'_j}M^i_aN^j_b$.

\bb{Quadratic algebras associated with idempotent operators.} Let $\B=TW/(R)\in\QSA$. The super-subspace $R\subset W\otimes W$ has a linear complement $R_c$, i.e. a super-subspace $R_c\subset W\otimes W$ such that $W\otimes W=R\oplus R_c$. This decomposition defines the operator $B\in\End(W\otimes W)$ by the formula $B(r+r_c)=r$ for $r\in R$, $r_c\in R_c$. It is an idempotent: $B^2=B$. Thus we have $\B=TW/(\Im B)$. Since $R_c$ is not unique for a fixed algebra $\B\in\QSA$, the idempotent $B$ is also not unique.

The entries $B^{st}_{ij}$ of the idempotent $B$ in the basis $(e_i\otimes e_j)$ are defined by the formula $B(e_i\otimes e_j)=\sum_{s,t\in I}B^{st}_{ij}(e_s\otimes e_t)=\sum_{s,t\in I}(e_s\otimes e_t)B^{st}_{ij}$. Since the operator $B\in\End(W\otimes W)$ is even, $B^{st}_{ij}=0$ if $\bar k_s+\bar k_t+\bar k_i+\bar k_j\ne\bar0$. The operator $B^*\in\End(W^*\otimes W^*)$ is also an idempotent. The dual basis is $(-1)^{\bar k_i\bar k_j}(e^i\otimes e^j)$, hence we have $B^*(e^s\otimes e^t)=(-1)^{\bar k_i\bar k_j+\bar k_s\bar k_t}B^{st}_{ij}(e^i\otimes e^j)$.

Any quadratic super-algebra $\B\in\FQSA$ has the form $\B\cong TW^*/(\Im B^*)$ for an idempotent operator $B\in\End(W\otimes W)$. We use the notations $\gX_B(\KK)=TW^*/(\Im B^*)$ and, more generally,
 $\gX_B(\R)=\R\otimes\gX_B(\KK)\cong\gX_B(\KK)\otimes\R\in\GrSAlg{\NN_0}$ for $\R\in\SAlg$.
 Denote by $x^i$ the element of $\gX_B(\KK)$ corresponding to the basis element $e^i\in W^*$. Then the super-algebra $\gX_B(\R)$ is an $\R$-algebra generated by $x^i$, $i\in I$, with parities $[x^i]=\bar k_i$ and commutation relations $\sum_{i,j\in I}(-1)^{\bar k_i\bar k_j}B^{st}_{ij}x^ix^j=0$.

Let $S=1-B$ be the dual idempotent: $\Im S=R_c$. Denote $\Xi_B(\KK)=TW/(\Im S)$, $\Xi_B(\R)=\Xi_B(\KK)\otimes\R\cong\R\otimes\Xi_B(\KK)$. The latter is an $\R$-algebra generated by $\psi_i$, $i\in I$, with parities $[\psi_i]=\bar k_i$ and commutation relations $\sum_{i,j\in I}S^{st}_{ij}\psi_s\psi_t=0$, where $\psi_i\in\Xi_B(\KK)$ corresponds to $e_i\in W$ and $S^{st}_{ij}=\delta^s_i\delta^t_j-B^{st}_{ij}$ are entries of $S$.

Consider the right $\gX_B(\R)$-module $W\otimes\gX_B(\R)$ and denote $X=\sum_{i\in I}e_ix^i\in W\otimes\gX_B(\R)$. We can consider $X$ as the column-vector $X=\begin{pmatrix}x^1 \\ \vdots \\ x^N\end{pmatrix}$. The element $\Psi=\sum_{i\in I}\psi_ie^i\in\Xi_B(\R)\otimes W^*$ of the left free $\Xi_B(\R)$-module $\Xi_B(\R)\otimes W^*$ can be regarded as the row-vector $\Psi=\begin{pmatrix}\psi_1 & \ldots & \psi_N\end{pmatrix}$. By considering $X$ as an operator $\KK\to W\otimes\gX_B(\R)$ we can use the notations~\eqref{M1}, \eqref{M2}, \eqref{M1N2}. We have $X^{(1)}X^{(2)}=(\id\otimes\id\otimes\mu_{\gX_B(\R)})(\sigma^{(23)})X\otimes X=\sum_{i,j\in I}(-1)^{\bar k_i\bar k_j}(e_i\otimes e_j)x^ix^j\in W\otimes W\otimes\gX_B(\R)$. Since $\Psi$ is an even element of $\Xi_B(\R)\otimes W^*\cong\big(W\otimes\Xi_B(\R)\big)^*$ we can regard it as an operator $W\to\Xi_B(\R)$, namely $\Psi(e_i)=\psi_i$, so $\Psi^{(1)}\Psi^{(2)}(e_i\otimes e_j)=\psi_i\psi_j$. The latter implies $\Psi^{(1)}\Psi^{(2)}=\sum_{i,j\in I}(-1)^{\bar k_i\bar k_j}\psi_i\psi_j(e^i\otimes e^j)\in\Xi_B(\R)\otimes W^*\otimes W^*$.
 The operator $B$ acts on $X^{(1)}X^{(2)}$ and $\Psi^{(1)}\Psi^{(2)}$ from the left and from the right respectively: $B(X^{(1)}X^{(2)})=\sum_{i,j\in I}(-1)^{\bar k_i\bar k_j}B(e_i\otimes e_j)x^ix^j=\sum_{i,j,s,t\in I}(-1)^{\bar k_i\bar k_j}B^{st}_{ij}(e_s\otimes e_t)x^ix^j$, $(\Psi^{(1)}\Psi^{(2)})B=\sum_{i,j,s,t\in I}(-1)^{\bar k_s\bar k_t}\psi_i\psi_jB^{ij}_{st}(e^s\otimes e^t)$. We see that the commutation relations of the $\R$-algebras $\gX_B(\R)$ and $\Xi_B(\R)$ are written in matrix form as
\begin{align}
 &B(X^{(1)}X^{(2)})=0, && (\Psi^{(1)}\Psi^{(2)})S=0.
\end{align}

Note that the quadratic super-algebras $\gX_B(\KK)$ and $\Xi_B(\KK)$ are Koszul-dual to each other: $\gX_B(\KK)^!=\Xi_B(\KK)$, $\Xi_B(\KK)^!=\gX_B(\KK)$.

Let $A_W=\frac{1-\sigma_{W,W}}2$, $S_W=\frac{1+\sigma_{W,W}}2$ be super-antisymmetrizer and super-symmetrizer acting in $W\otimes W$. Then we have
\begin{align}
 &SW^*=\gX_{A_W}(\KK), &&\Lambda W=\Xi_{A_W}(\KK).
\end{align}

\begin{Lem} \label{LemTXX}
 Let $T\colon W\otimes W\to\wt W\otimes\R$, then the condition $T(X^{(1)}X^{(2)})=0$ is equivalent to $TS=0$. Let $U\colon\wt W\to W\otimes W\otimes\R$, then the condition $(\Psi^{(1)}\Psi^{(2)})U=0$ is equivalent to $BU=0$.
\end{Lem}

\noindent{\bf Proof.} Define the entries $T^a_{ij}\in\R$ as $T(e_i\otimes e_j)=\sum_{a\in\wt I}\wt e_aT^a_{ij}$. If $TS=0$, then $T=TB$ and hence $T(X^{(1)}X^{(2)})=TB(X^{(1)}X^{(2)})=0$. Conversely, let the element $T(X^{(1)}X^{(2)})=\sum_{i,j\in I}(-1)^{\bar k_i\bar k_j}T(e_i\otimes e_j)x^ix^j=\sum_{i,j\in I\atop a\in\wt I}(-1)^{\bar k_i\bar k_j}\wt e_aT^a_{ij}x^ix^j$ vanish, then $\sum_{i,j\in I}(-1)^{\bar k_i\bar k_j}T^a_{ij}x^ix^j=0$, hence $T^a_{ij}=G^a_{st}B^{st}_{ij}$ for some $G^a_{st}\in\R$. We get $T=GB$, where $G\colon W\otimes W\to\R\otimes\wt W$ is the operator with the entries $G^a_{ij}$. Thus $TS=GBS=0$. The second equivalence is checked similarly. \qed

\bb{$(B,\wt B)$-Manin matrices.}
Let $\wt B\in\End(\wt W\otimes\wt W)$ be another idempotent and $\wt S=1-\wt B$. Their entries are defined by $\wt B(\wt e_c\otimes\wt e_d)=\sum_{a,b\in\wt I}(\wt e_a\otimes\wt e_b)\wt B^{ab}_{cd}$, $\wt S^{ab}_{cd}=\delta^a_c\delta^b_d-\wt B^{ab}_{cd}$. Denote the generators of $\gX_{\wt B}(\KK)$ and $\Xi_{\wt B}(\KK)$ by $\wt x^a$ and $\wt\psi_a$, where $a$ runs over $\wt I$. They form the column-vector $\wt X=\sum_{a\in\wt I}\wt e_a \wt x^a$ and the row-vector $\wt\Psi=\sum_{a\in\wt I}\wt\psi_a \wt e^a$.

Any morphism of $\NN_0$-graded super-algebras $f\colon\gX_B(\KK)\to\gX_{\wt B}(\R)$ is uniquely determined by its first component $f_1\colon W^*\to\R\otimes\wt W^*$. In terms of the generators we have $f(x^i)=\sum_{a\in\wt I}M^i_a \wt x^a$ for some matrix $M=(M^i_a)$. This is a matrix of the operator $M\colon\wt W\to\R\otimes W$ such that $f_1=M^\st$. Analogously, a morphism $g\colon\Xi_{\wt B}(\KK)\to\Xi_B(\R)$ in $\GrSAlg{\NN_0}$ has the form $g(\wt\psi_a)=\sum_{i\in I}\psi_iM^i_a$.

\begin{Th} \label{ThAAManin}
 The following conditions are equivalent.
\begin{itemize}
 \item The operator $M\colon\wt W\to W\otimes\R$ satisfies
\begin{align} \label{AMMA}
 BM^{(1)}M^{(2)}(1-\wt B)=0.
\end{align}
 \item The entries $M^i_a\in\R$ have the parities $[M^i_a]=\bar k_i+\bar l_a$ and satisfy
\begin{align} \label{AMMAij}
 \sum_{i,j\in I\atop a,b\in\wt I}(-1)^{(\bar k_i+\bar l_a)\bar k_j}B^{st}_{ij}M^i_aM^j_b\wt S^{ab}_{cd}=0.
\end{align}
 \item The formula
\begin{align} 
 f_M(x^i)=\sum_{a\in\wt I}M^i_a\wt x^a
\end{align}
defines a graded morphism $f_M\colon\gX_B(\KK)\to\gX_{\wt B}(\R)$.
 \item The formula
\begin{align} 
 f^M(\wt\psi_a)=\sum_{i\in I}\psi_iM^i_a
\end{align}
defines a graded morphism $f^M\colon\Xi_{\wt B}(\KK)\to\Xi_{B}(\R)$.
\end{itemize}
\end{Th}

\noindent{\bf Proof.} The first and second conditions are equivalent since the left hand side of~\eqref{AMMAij} is the entries of the left hand side of~\eqref{AMMA}. Let $y^i=\sum_{a\in\wt I}M^i_a\wt x^a$ and $Y=\sum_{i\in I}e_i y^i=M\wt X$.
 The formula $f_M(x^i)=y^i$ defines the homomorphism iff $\sum_{i,j\in I}(-1)^{\bar k_i\bar k_j}B^{st}_{ij}y^iy^j=0$, that is $B(Y^{(1)}Y^{(2)})=0$. Substitution $Y=M\wt X$ yields $B(M^{(1)}M^{(2)})(\wt X^{(1)}\wt X^{(2)})=0$. By means of Lemma~\ref{LemTXX} applied to $T=B(M^{(1)}M^{(2)})$ this is equivalent to~\eqref{AMMA}. Similarly, the last condition means exactly that $(\Phi^{(1)}\Phi^{(2)})\wt S=0$, where $\Phi=\Psi M$. By virtue of Lemma~\ref{LemTXX} applied to $U=(M^{(1)}M^{(2)})\wt S$ this is also equivalent to~\eqref{AMMA}. \qed

\begin{Def} \normalfont
 The matrix $M=(M^i_a)$ satisfying the conditions from Theorem~\ref{ThAAManin} is called {\it super-Manin matrix for idempotents} $B$ and $\wt B$ or {\it $(B,\wt B)$-Manin matrix}. If $B=\wt B$, then we call it simply {\it $B$-Manin matrix}.
\end{Def}

Thus we obtain a one-to-one correspondence between $(B,\wt B)$-Manin matrices, the morphisms $f_M\colon\gX_B(\KK)\to\gX_{\wt B}(\R)$ and the morphisms $f^M\colon\Xi_{\wt B}(\KK)\to\Xi_{B}(\R)$. These morphisms can be written in the matrix form:
\begin{align}
 &f_M(X)=M\wt X, && f^M(\wt\Psi)=\Psi M.
\end{align}

\bb{Universal $(B,\wt B)$-Manin matrices.} \label{bbUMM}
Let $\U_{B,\wt B}$ be the quadratic super-algebra with the generators $\M^i_a$, $i\in I$, $a\in\wt I$, of parity $[\M^i_a]=\bar k_i+\bar l_a$; the commutation relations are $\sum_{i,j\in I\atop a,b\in\wt I}(-1)^{(\bar k_i+\bar l_a)\bar k_j}B^{st}_{ij}\M^i_a\M^j_b\wt S^{ab}_{cd}=0$. These generators form the $(B,\wt B)$-Manin matrix $\M=(\M^i_a)$. We call $\M$ the {\it universal $(B,\wt B)$-Manin matrix} by the following reason: any $(B,\wt B)$-Manin matrix $M$ over $\R$ is an image of $\M$ under a homomorphism of super-algebras $h\colon\U_{B,\wt B}\to\R$ in the sense that $M^i_a=h(\M^i_a)$. Theorem~\ref{ThAAManin} gives the bijections
\begin{align} \label{bijUAA}
 \Hom_\SAlg(\U_{B,\wt B},\R)=\Hom_\FQSA\big(\gX_B(\KK),\gX_{\wt B}(\R)\big)=\Hom_\FQSA\big(\Xi_{\wt B}(\KK),\Xi_B(\R)\big).
\end{align}

\begin{Prop} \label{PropUMM}
We have the isomorphisms of quadratic super-algebras
\begin{align} \label{Ucohom}
\U_{B,\wt B}\cong\bfcohom\big(\gX_{\wt B}(\KK),\gX_B(\KK)\big)\cong\bfcohom\big(\Xi_B(\KK),\Xi_{\wt B}(\KK)\big)
\end{align}
given by the identification $\M^i_a=x^i\otimes\wt\psi_a=(-1)^{\bar k_i\bar l_a}\wt\psi_a\otimes x^i$.
\end{Prop}

\noindent{\bf Proof.} Recall that $\bfcohom\big(\gX_B(\KK),\gX_{\wt B}(\KK)\big)$ has the form $\gX_B(\KK)\bullet\gX_{\wt B}(\KK)^!=\gX_B(\KK)\bullet\Xi_{\wt B}(\KK)$. Let $M^i_a=x^i\otimes\wt\psi_a$. If $\sum\limits_{i,j,a,b}T^{ab}_{ij}M^i_aM^j_b=0$ for some $T^{ab}_{ij}\in\KK$, then $\sum\limits_{i,j,a,b}T^{ab}_{ij}(-1)^{\bar k_j\bar l_a}(e^i\otimes e^j\otimes\wt e_a\otimes\wt e_b)\in\Im(B^*)\otimes\Im(\wt S)$, so we have $T^{ab}_{ij}=\sum\limits_{s,t,c,d}G^{cd}_{st}(-1)^{\bar k_j(\bar l_a+\bar k_i)}B^{st}_{ij}\wt S^{ab}_{cd}$ for some $G^{cd}_{st}\in\KK$, hence the relation $\sum\limits_{i,j,a,b}T^{ab}_{ij}M^i_aM^j_b=0$ is a linear combination of the relations~\eqref{AMMAij}. This implies the first isomorphism~\eqref{Ucohom}. Then, the second isomorphism~\eqref{Ucohom} follows from~\eqref{BAAB}. \qed

Let us generalise the result~\cite[\S~4.2.3, Th.~4.3]{Sgrt} to the super-case by using the quadratic super-algebra $\U_{B,\wt B}$.

\begin{Th} \label{ThCoclosed}
 The symmetric monoidal categories $(\GrSAlg{\NN_0},\circ)$ and $(\GrSAlg{\ZZ},\circ)$ are coclosed relative to $\bfP=\FQSA$. The cohom-functor $\bfcohom\colon\bfP^\op\times\bfP\to\bfC$ is the composition of the cohom-functor $\bfcohom\colon\bfP^\op\times\bfP\to\QSA$ described in Propositions~\ref{Propcohom}, \ref{PropUMM} with the embedding $\QSA\hookrightarrow\bfC$ for the both $\bfC=\GrSAlg{\NN_0}$ and $\bfC=\GrSAlg{\ZZ}$.
\end{Th}

\noindent{\bf Proof.} If $\R$ has a structure of $\ZZ$-graded super-algebra $\R=\bigoplus\limits_{k\in\ZZ}\R_k$, then the graded homomorphisms $h\colon\U_{B,\wt B}\to\R$ form a subset of~\eqref{bijUAA} corresponding to the $(B,\wt B)$-Manin matrices $M$ such that $M^i_a\in\R_1$. As a subset of $\Hom_\FQSA\big(\gX_B(\KK),\gX_{\wt B}(\R)\big)$ it consists of the homomorphisms factorising as $f_M\colon\gX_B(\KK)\xrightarrow{f}\R\circ\gX_{\wt B}(\KK)\hookrightarrow\gX_{\wt B}(\R)$ through the graded homomorphism $f(x^i)=\sum\limits_{a\in\wt I}M^i_a\otimes\wt x^a$.
 Consider the graded homomorphism $\eta\colon\gX_B(\KK)\to\U_{B,\wt B}\circ\gX_{\wt B}(\KK)$, $\eta(x^i)=\sum\limits_{a\in\wt I}\M^i_a\otimes\wt x^a$. For any morphism $f$ there exists a unique $h$ making the diagram
\begin{align}
\xymatrix{
 \gX_B(\KK)\ar[rr]^{\eta\qquad}\ar[rrd]_f && \U_{B,\wt B}\circ\gX_{\wt B}(\KK)\ar[d]^{h\circ\id} \\
 &&\R\circ\gX_{\wt B}(\KK) 
}
\end{align}
commute. Hence $(\U_{B,\wt B},\eta)$ is a universal morphism from the object $\gX_B(\KK)$ to the functor $G\colon\GrSAlg{\ZZ}\to\GrSAlg{\ZZ}$, $G(\R)=\R\circ\gX_{\wt B}(\KK)$. Since any $\B,\wt\B\in\FQSA$ are isomorphic to some $\gX_B(\KK)$ and $\gX_{\wt B}(\KK)$, the monoidal category $(\GrSAlg{\ZZ},\circ)$ is coclosed relative to $\FQSA$ with the cohom-functor~\eqref{Ucohom}, see~\cite[\S~2.3.2, Th.~3.13 (2')]{Sgrt}. By restricting this proof to the case $\R=\bigoplus\limits_{k\ge0}\R_k$ we obtain the same result for $(\GrSAlg{\NN_0},\circ)$. \qed

Consider the case $\B=\wt\B=\gX_B(\KK)$, then $\M=(\M^i_j)_{i,j\in I}$ is a universal $B$-Manin matrix. The quadratic algebra $\bfcoend(\B)=\U_B:=\U_{B,B}$ has a structure of comonoid in $(\GrSAlg{\NN_0},\circ)$ given by the graded homomorphisms
\begin{align}
 &d_\B\colon\bfcoend(\B)\to\bfcoend(\B)\circ\bfcoend(\B), &&\M^i_j\mapsto\sum_{l\in I}\M^i_l\otimes\M^l_j, \label{dB} \\
 &v_\B\colon\bfcoend(\B)\to\KK[u], &&\M^i_j\mapsto\delta^i_j u. \label{vB}
\end{align}

\subsection{Quantum representations and quantum linear actions}
\label{sec42}

\bb{(Co)representations.}
Let $(\bfC,\otimes)$ be a relatively closed monoidal category with the parametrising subcategory $\bfP\subset\bfC$. By following~\cite{Sqrt}, \cite{Sgrt} we define {\it corepresentation} of a comonoid $\OO=(X,\Delta_X,\varepsilon_X)\in\Comon(\bfC,\otimes)$ on an object $W\in\bfP$ as a morphism $\omega\colon\bfcoend(W)\to\OO$ in the category $\Comon(\bfC,\otimes)$. Corepresentations of a fixed $\OO$ form the category $\Corep_\bfP(\OO)$. Its objects are $(W,\omega)$. A morphism $(W,\omega)\to(W',\omega')$ in $\Corep_\bfP(\OO)$ is a morphism $f\colon W\to W'$ in $\bfC$ such that the diagram
\begin{align} \label{omegaf}
\xymatrix{
 \bfcohom(W',W)\ar[rrr]^{\bfcohom(f,\id_W)}\ar[d]^{\bfcohom(\id_{W'},f)} &&& \bfcohom(W,W)\ar[d]^{\omega} \\
\bfcohom(W',W')\ar[rrr]^{\qquad\omega'} &&& X
}
\end{align}
commutes. The adjunction morphism
\begin{align} \label{varth}
 \vartheta\colon\Hom\big(\bfcoend(W),X\big)\isoright\Hom(W,X\otimes W)
\end{align}
gives a one-to-one correspondence between the corepresentations $\omega\colon\bfcoend(W)\to\OO$ and coactions $\delta\colon W\to X\otimes W$ of $\OO$ on $W$. In this way we obtain the embedding of categories $\Corep_\bfP(\OO)\hookrightarrow\Lcoact(\OO)$.

Dually, if $(\bfC,\otimes)$ is closed relative to $\bfP$, then one can define representations of a monoid $\MM\in\Mon(\bfC,\otimes)$. They form the category $\Rep_\bfP(\MM)$ consisting of $(W,\rho)$, where $W\in\bfP$ and $\rho\colon\MM\to\bfend(W)$.
For example, if $(\bfC,\otimes)$ is the closed symmetric monoidal category $(\SVect,\otimes)$ defined in p.~\ref{bbSVect} (the case $\bfP=\bfC=\SVect$), then we obtain the categories $\Rep_\SVect(\R)$ for $\R\in\SAlg=\Mon(\SVect,\otimes)$. This is a category of left super-modules over the super-algebra $\R$.

\bb{Quantum super-algebras.}
Let us consider the cases of the relatively closed categories $(\GrSAlg{\NN_0},\circ)$ and $(\GrSAlg{\ZZ},\circ)$. A comonoid in $(\GrSAlg{\ZZ},\circ)$ is $\OO=(\A,\Delta,\varepsilon)$, where $\A\in\GrSAlg{\ZZ}$ is a $\ZZ$-graded super-algebra $\A=\bigoplus\limits_{k\in\ZZ}\A_k$ and $\Delta\colon\A\to\A\circ\A$ and $\varepsilon\colon\A\to\KK[u,u^{-1}]$ are graded homomorphisms with components $\Delta_k\colon\A_k\to\A_k\otimes\A_k$, $\varepsilon_k\colon\A_k\to\KK u^k$. By composing them with non-graded homomorphisms $\A\circ\A\hookrightarrow\A\otimes\A$ and $\KK[u,u^{-1}]\to\KK$, $f(u)\mapsto f(1)$ we obtain morphisms $\Delta_\A\colon\A\to\A\otimes\A$ and $\varepsilon_\A\colon\A\to\KK$ in $\SAlg$ such that $(\A,\Delta_\A,\varepsilon_\A)$ is a comonoid in $(\SAlg,\otimes)$, i.e. a super-bialgebra. Conversely, a super-bialgebra $(\A,\Delta_\A,\varepsilon_\A)$ gives a comonoid $\OO=(\A,\Delta,\varepsilon)$ if $\A$ is a $\ZZ$-graded super-algebra satisfying
\begin{align} \label{DeltaAk}
 \Delta_\A(\A_k)\subset\A_k\otimes\A_k;
\end{align}
the morphism $\varepsilon\colon\A\to\KK[u,u^{-1}]$ has the form $\varepsilon(a_k)=\varepsilon_\A(a_k)u^k$, where $a_k\in\A_k$.

Analogously, a comonoid $\OO=(\A,\Delta,\varepsilon)\in\Comon(\GrSAlg{\NN_0},\circ)$ is given by $\A=\bigoplus\limits_{k\ge0}\A_k$, $\Delta\colon\A\to\A\circ\A$, $\varepsilon\colon\A\to\KK[u]$. Such comonoids correspond to bialgebras $(\A,\Delta_\A,\varepsilon_\A)$ with $\NN_0$-grading of $\A$ making it an $\NN_0$-graded super-algebra satisfying~\eqref{DeltaAk}. The colax monoidal embedding $(\GrSAlg{\NN_0},\circ)\hookrightarrow(\GrSAlg{\ZZ},\circ)$ induces the categorical embedding $\Comon(\GrSAlg{\NN_0},\circ)\subset\Comon(\GrSAlg{\ZZ},\circ)$. The comonoids in $(\QSA,\circ)$ (more precisely, monoids in $(\QSA^\op,\circ)$) generalise the super-algebras for the quantum case. By {\it quantum super-algebra} we mean a general comonoid $\OO\in\Comon(\GrSAlg{\ZZ},\circ)$.

The comonoid $\bfcoend(\B)\in\Comon(\GrSAlg{\NN_0},\circ)$ described in p.~\ref{bbUMM} is embedded into $\Comon(\GrSAlg{\ZZ},\circ)$ as the $\ZZ$-graded super-algebra $\bfcoend(\B)=\U_{B}$ with the comultiplication~\eqref{dB} and the counit $v_\B\colon\bfcoend(\B)\to\KK[u,u^{-1}]$ defined by the formula~\eqref{vB}.

\bb{Quantum representations and actions.}
Let $\OO=(\A,\Delta,\varepsilon)$ be a quantum super-algebra.

\begin{Def} \normalfont
 {\it Quantum representation} of $\OO$ on $\B\in\FQSA$ is the corepresentation $\omega\colon\bfcoend(\B)\to\OO$. This is a graded homomorphism $\omega\colon\bfcoend(\B)\to\A$ making the diagrams
\begin{align}
 &\xymatrix{
\bfcoend(\B)\ar[d]^\omega\ar[r]^{d_\B\quad\qquad} & \bfcoend(\B)\circ\bfcoend(\B)\ar[d]^{\omega\circ\omega} \\
\A\ar[r]^\Delta & \A\circ\A
} &&
\xymatrix{
\bfcoend(\B)\ar[d]^\omega\ar[rd]^{v_\B} \\
\A\ar[r]^{\varepsilon\qquad} & \KK[u,u^{-1}]
}
\end{align}
commute. Their morphisms are morphisms in $\Corep_\FQSA(\OO)$.
 {\it Quantum linear action} of $\OO$ on $\C\in\GrSAlg{\ZZ}$ is the coaction $\delta\colon\C\to\A\circ\C$. This is a graded homomorphism such that the following diagrams commute:
\begin{align}
 &\xymatrix{
\C\ar[d]^\delta\ar[rr]^{\delta} && \A\circ\C\ar[d]^{\Delta\circ\id} \\
\A\circ\C\ar[rr]^{\id\circ\delta} && \A\circ\A\circ\C
} &&
\xymatrix{
\C\ar[d]^\delta\ar@{=}[rrd] \\
\A\circ\C\ar[rr]^{\varepsilon\circ\id\qquad} && \KK[u,u^{-1}]\circ\C
}
\end{align}
Their morphisms are morphisms in $\Lcoact(\OO)$.
\end{Def}

For a quantum linear action $\delta$ consider a composition $\delta_\A\colon\C\xrightarrow{\delta}\A\circ\C\hookrightarrow\A\otimes\C$. This is a coaction of the super-bialgebra $(\A,\Delta_\A,\varepsilon_\A)$ on the super-algebra $\C$. Conversely, a coaction $\delta_\A\colon\C\to\A\otimes\C$ of this super-bialgebra has this form for some $\delta$ if
\begin{align}
 \delta_A(\C_k)\subset\A_k\otimes\C_k.
\end{align}
Thus we obtain the forgetful functor $\Lcoact(\OO)\to\Lcoact(\A,\Delta_\A,\varepsilon_\A)$.

The bijection~\eqref{varth} for $(\GrSAlg{\ZZ},\circ)$ has the form
\begin{align} \label{varthZZ}
 \vartheta\colon\Hom\big(\bfcoend(\B),\A\big)\isoright\Hom(\B,\A\circ\B).
\end{align}
For any quantum super-algebra $\OO=(\A,\Delta,\varepsilon)$ and quadratic super-algebra $\B\in\FQSA$ we have the one-to-one correspondence $\vartheta\colon\omega\leftrightarrow\delta$ between quantum representations $\omega$ and quantum linear actions $\delta$ of $\OO$ on $\B$. It defines the fully faithful functor $\Corep_\FQSA(\OO)\hookrightarrow\Lcoact(\OO)$, $(\B,\omega)\mapsto\big(\B,\vartheta(\omega)\big)$.

\bb{Multiplicative Manin matrices.}
The matrix $M=(M^i_j)$ over a super-bialgebra $(\A,\Delta_\A,\varepsilon_\A)$ is called {\it multiplicative} iff
\begin{align}
 &\Delta_\A(M^i_j)=\sum_l M^i_l\otimes M^l_j, && \varepsilon_\A(M^i_j)=\delta^i_j.
\end{align}
Let $\OO=(\A,\Delta,\varepsilon)$ be a quantum super-algebra and $(\A,\Delta_\A,\varepsilon_\A)$ be the corresponding super-bialgebra. The matrix $M$ with entries $M^i_j\in\A$ is called multiplicative over $\OO$ iff it is multiplicative over this super-bialgebra. It is called {\it first order} matrix iff $M^i_j\in\A_1$. For instance, the universal $B$-Manin matrix is a first order multiplicative matrix over $\U_B=\bfend\big(\gX_B(\KK)\big)$.

\begin{Th}
 Let $\OO\in\Comon(\GrSAlg{\ZZ},\circ)$ and $\B=\gX_B(\KK)$ for some idempotent $B$. Let $\M=(\M^i_j)$ be the universal $B$-Manin matrix. Any quantum representation of $\OO$ on $\B$ has the form $\omega(\M^i_j)=M^i_j$ for some first order multiplicative $B$-Manin matrix $M=(M^i_j)$ over $\OO$. Conversely, any such matrix defines a quantum representation $\omega\colon\bfcoend(\B)\to\OO$. The quantum linear action $\delta=\vartheta(\omega)\colon\B\to\A\circ\B$ has the form $\delta(x^i)=\sum_j M^i_j\otimes x^j$, i.e. $\delta_\A=f_M$. Let $\C=\gX_C(\KK)$. Consider two quantum representations $\omega\colon\bfcoend(\B)\to\OO$ and $\nu\colon\bfcoend(\C)\to\OO$ corresponding to the multiplicative first order $B$- and $C$-Manin matrices $M$ and $N$ over $\OO$. A homomorphism $f\colon\B\to\C$ is a morphism $(\B,\omega)\to(\C,\nu)$ in $\Corep_\FQSA(\OO)$ (=morphism $\big(\B,\vartheta(\omega)\big)\to\big(\C,\vartheta(\nu)\big)$ in $\Lcoact(\OO)$) iff $f=f_K$ for a $(B,C)$-Manin matrix $K$ over $\KK$ such that $MK=KN$.
\end{Th}

\noindent This theorem follows from Theorems~\ref{ThAAManin} and~\ref{ThCoclosed}. The proof is similar to~\cite[Th.~5.2, Prop.~5.3,~5.5,~5.6]{Sqrt}.

\bb{Opposite and coopposite quantum representations.}
Let $\R=(V,\mu_\R,\eta_\R)$ be a super-algebra with a multiplication $\mu_\R\colon V\otimes V\to V$ and a unity map $\eta_\R\colon V\to\KK$, where $V\in\SVect$. The opposite super-algebra is $\R^\op=(V,\mu_\R^\op,\eta_\R)$. If $\R\in\GrSAlg{\ZZ}$, then $\R^\op\in\GrSAlg{\ZZ}$, see p.~\ref{bbFunctop}. 

The opposite and coopposite to the comonoid $\OO=(\A,\Delta,\varepsilon)\in\Comon(\GrSAlg{\ZZ},\circ)$ are $\OO^\op=(\A^\op,\Delta,\varepsilon)\in\Comon(\GrSAlg{\ZZ},\circ)$, $\OO^\cop=(\A,\Delta^\cop,\varepsilon)\in\Comon(\GrSAlg{\ZZ},\circ)$, where $\Delta^\cop=\sigma\cdot\Delta$. The bialgebras corresponding to these comonoids are $(\A^\op,\Delta_\A,\varepsilon_\A)$ and $(\A,\Delta_\A^\cop,\varepsilon_\A)$.

Denote $B^{(21)}=\sigma\cdot B\sigma$. This is an operator $W\otimes W\to W\otimes W$ with the entries $(B^{(21)})^{st}_{ij}=(-1)^{\bar k_s\bar k_t+\bar k_i\bar k_j}B^{ts}_{ji}$. The quadratic super-algebra $\gX_{B^{(21)}}(\KK)$ is defined by the commutation relations $\sum_{i,j}B^{ts}_{ji}x^ix^j=0$, hence $\gX_{B^{(21)}}(\KK)=\gX_B(\KK)^\op$.
Also, we obtain $\U_{B^{(21)}}=\U_B^\op$, which is an isomorphism of quantum super-algebras: $\bfcoend(\B^\op)=\bfcoend(\B)^\op$.

Consider a quantum representation $\omega\colon\bfcoend(\B)\to\OO$ of $\OO=(\A,\Delta,\varepsilon)$ on $\B=\gX_B(\KK)$ corresponding to a multiplicative first order $B$-Manin matrix $M$ over $\A$.
 It gives the comonoid morphism $\omega^\op\colon\bfcoend(\B^\op)\to\OO^\op$. This is a quantum representation of $\OO^\op$ on $\B^\op=\gX_{B^{(21)}}(\KK)$ defined by the multiplicative first order $B^{(21)}$-Manin matrix $M$ over $\A^\op$.

Due to~\eqref{BAAB} the quantum super-algebra $\bfcoend(\B)$ coincides with $\bfcoend(\B^!)$ as a graded super-algebra, but it has different comultiplication $d_{\B^!}$. Note that the matrix $\M^\ist$ is multiplicative over the comonoid $\bfcoend(\B^!)$, i.e. $d_{\B^!}\colon(M^\ist)^i_j\mapsto\sum_l(M^\ist)^i_l\otimes(M^\ist)^l_j$, where $(M^\ist)^i_j=(-1)^{(\bar k_j+\bar k_i)\bar k_j}M^j_i$ are entries of $\M^\ist$, see p.~\ref{bbMatrOp}. Hence we have $\bfcoend(\B^!)=\bfcoend(\B)^\cop$.
The quantum representation $\omega\colon\bfcoend(\B)\to\OO$ gives the comonoid morphism $\omega^\cop\colon\bfcoend(\B^!)\to\OO^\cop$. This is a quantum representation of $\OO^\cop$ on $\B^!=\Xi_{B}(\KK)$ defined by the multiplicative first order $B^*$-Manin matrix $M^\ist$ over $\A^\cop$. The corresponding coaction is $\delta^\cop(\psi_i)=\sum_j(M^\ist)^i_j\psi_j=\sum_j\psi_jM^j_i$.

\bb{Parity change.}
Let us apply the functor $\wh\Pi\colon\FQSA\to\FQSA$ to the quantum representation space $\B=\gX_B(\KK)$. By definition we need to change $W$ by $\Pi W$, which is equivalent to the substitution $\bar k_i\to\bar k_i+\bu$. Denote by $\Pi B$ the operator $\Pi W\otimes\Pi W\to\Pi W\otimes\Pi W$ given by the same matrix, but with the changed format. Since $B^{st}_{ij}\ne0$ implies $\bar k_s+\bar k_t=\bar k_i+\bar k_j$ the relations $B^{st}_{ij}(-1)^{\bar k_1\bar k_j}x^ix^j=0$ keep the same, so that $\wh\Pi\gX_B(\KK)=\gX_{\Pi B}(\KK)$.

We do not change the $\ZZ_2$-grading of $\A$. If $M=(M^i_j)$ is a multiplicative first order $B$-Manin matrix over $\A$, then the parity $[M^i_j]=\bar k_i+\bar k_j$ does not change, hence $M$ is a $\Pi B$-Manin matrix. Thus we obtain a quantum representation of the same $\OO=(\A,\Delta,\varepsilon)$ on the quadratic super-algebra $\wh\Pi\B=\gX_{\Pi}(\KK)$ defined by the multiplicative first order $\Pi B$-Manin matrix $M$.

\bb{Classical left modules over finite-dimensional super-algebras.}
Let $\R$ be a finite-dimensional super-algebra, that is $\R\in\Mon(\FSVect,\otimes)$. It can be lifted to the level of quantum super-algebras by the contravariant colax monoidal functor $S^*\colon(\FSVect,\otimes)\to(\FQSA,\circ)$. We obtain the quantum algebra $S^*\R=S\R^*$. This is a commutative super-algebra of polynomials on the super-vector space $\R$. The multiplication and unity of $\R$ induces the structure of super-bialgebra on $S\R^*$.

A structure of (left) $\R$-module on a super-vector space $W$ is given by a super-algebra morphism $\rho\colon\R\to\bfend(W)$ or by an action $a\colon\R\otimes W\to W$. In the basis $(e_i)$ they have the form $\rho(r)e_j= a(r\otimes e_j)=\sum_i\rho^i_j(r)e_i$ for some linear functions $\rho^i_j\in\R^*$ such that $[\rho^i_j]=\bar k_i+\bar k_j$ and the matrix $M=(\rho^i_j)$ is multiplicative over $S\R^*$. Since $S\R^*$ is a commutative super-algebra this matrix is an $A_W$-Manin matrix. The functor $S^*$ translates the representation $\rho$ to the quantum representation $\omega\colon\bfcoend(SW^*)\to S\R^*$ of $S\R^*$ on $S^*W=SW^*=\gX_{A_W}(\KK)$ (see details in~\cite{Sqrt,Sgrt}). This quantum representation is defined by the multiplicative first order $A_W$-Manin matrix $M=(\rho^i_j)$. The functor $S^*$ translates the action $a$ to a coaction $\delta\colon SW^*\to S\R^*\circ SW^*$.

Alternatively one can also use the contravariant strong monoidal functor $T^*$. In this case the classical representation $\rho$ is lifted to the quantum representation of the quantum algebra $T^*\R=T\R^*$ on the quadratic super-algebra $T^*W=TW^*=\gX_0(\KK)$ given by the same matrix $M=(\rho^i_j)$.

\bb{Quantum representations of quantum super-monoids.} A quantum super-monoid (super-group) is a super-bialgebra (Hopf super-algebra) $\BB\in\Bimon(\SVect,\otimes)=\Comon(\SAlg,\otimes)$. The strong monoidal functor
\begin{align}
 &(\SAlg,\otimes)\to(\GrSAlg{\ZZ},\circ), &&\R\mapsto\R\otimes\KK[u,u^{-1}]=\bigoplus_{k\in\ZZ}\R
\end{align}
is symmetric and faithful. It induces the non-full categorical embedding $\Bimon(\SVect,\otimes)\to\Comon(\GrSAlg{\ZZ},\circ)$, which identifies $\BB=(\R,\Delta_\R,\varepsilon_\R)$ with the quantum algebra $\OO_\BB:=\big(\R\otimes\KK[u,u^{-1}],\Delta_\R\otimes\KK[u,u^{-1}],\varepsilon_\R\otimes\KK[u,u^{-1}]\big)$.

By {\it quantum representation of a quantum monoid} $\BB$ on $\B=\gX_B(\KK)$ we mean a quantum representation of the corresponding comonoid $\OO_\BB$, i.e. a graded homomorphism of bialgebras $\omega\colon\bfcoend(\B)\to\OO_\BB$. It is given by a multiplicative $B$-Manin matrix $M$ over $\BB$.

\section*{Concluding remarks}
\addcontentsline{toc}{section}{Concluding remarks}

We generalised Manin's theory of quadratic algebras and Quantum Representation Theory introduced in~\cite{Sqrt} to the super-case. To do the latter we applied the general approach described in~\cite{Sgrt} to the monoidal category $(\GrSAlg{\ZZ},\circ)$ with the parametrising subcategory $\FQSA$. The abstract algebraic formulation and the language of monoidal categories allowed us to write all the formulae in the exactly the same form as for the non-super case~\cite{Sqrt}. We expect that the approach~\cite{Sgrt} works well for more general cases. The list of planed developments of the theory is written in the end of the article~\cite{Sqrt}.

\end{document}